\input amstex

\documentstyle{amsppt}
\magnification=\magstep1

\hsize=5.2in
\vsize 7in
\topmatter
\centerline {\bf  UNIQUENESS OF THE GROUP MEASURE SPACE}
\vskip 0.05in
\centerline {\bf DECOMPOSITION FOR POPA'S $\Cal H\Cal T$ FACTORS}
\rightheadtext{Group measure space decompositions of $\Cal H\Cal T$ factors}

\vskip 0.15in
\centerline { ADRIAN IOANA \footnote{Supported by a Clay Research Fellowship.}}
\address Department of Mathematics, University of California, San Diego, La Jolla, CA
92093-0112, USA 
\endaddress
\email aioana\@ucsd.edu \endemail

\thanks 2010 {\it Mathematics Subject Classification.} Primary 46L36; Secondary 28D15, 37A20.
\endthanks

\abstract We prove that if $\Gamma\curvearrowright (X,\mu)$ is a free ergodic rigid (in the sense of [Po01]) probability measure preserving action 
of a group $\Gamma$ with positive first $\ell^2$--Betti number, then the II$_1$ factor $L^{\infty}(X)\rtimes\Gamma$ has a unique group measure space Cartan subalgebra, up to unitary conjugacy.  We deduce that many $\Cal H\Cal T$ factors, including the II$_1$ factors associated with the usual actions $\Gamma\curvearrowright \Bbb T^2$ and $\Gamma\curvearrowright$ SL$_2(\Bbb R)$/SL$_2(\Bbb Z)$, where $\Gamma$ is a non--amenable subgroup of SL$_2(\Bbb Z)$, have a unique group measure space decomposition.
\endabstract

\endtopmatter
\document

\head \S 0. {Introduction and statement of the main results}.\endhead
\vskip 0.2in

The {\it group measure space construction} associates to every probability measure preserving (p.m.p.) action $\Gamma\curvearrowright (X,\mu)$ of a countable group $\Gamma$, a finite von Neumann algebra $M=L^{\infty}(X)\rtimes\Gamma$ ([MvN36]). If the action is free and ergodic, then $M$  is a II$_1$ factor and $A=L^{\infty}(X)$ is a {\it Cartan subalgebra}, i.e. a maximal abelian von Neumann subalgebra whose normalizer, $\Cal N_{M}(A)=\{u\in\Cal U(M)|uAu^*=A\}$, generates $M$.

During the last decade, S. Popa's {\it deformation/rigidity} theory has led to spectacular progress in the study of  II$_1$ factors (see  the surveys [Po07],[Va10a]).  In particular, several large families of group measure space II$_1$ factors $L^{\infty}(X)\rtimes\Gamma$  have been shown to have a unique Cartan subalgebra ([OP07],[OP08],[CS11]) or group measure space Cartan subalgebra ([Pe09],[PV09],[Io10],[FV10],[IPV10],[CP10],[HPV10],[Va10b]), up to unitary conjugacy. 
Such ``unique Cartan subalgebra" results play a crucial role in the classification of group measure space factors. More precisely,
they allow one to reduce the classification of the factors $L^{\infty}(X)\rtimes\Gamma$,  up to isomorphism, to the classification of the corresponding actions $\Gamma\curvearrowright X$, up to {\it orbit equivalence}. 
Indeed, by [Si55],[FM77], an isomorphism of group measure space factors $L^{\infty}(X)\rtimes\Gamma\cong L^{\infty}(Y)\rtimes\Lambda$ which identifies the Cartan subalgebras $L^{\infty}(X),L^{\infty}(Y)$, must come from an orbit equivalence between the actions, i.e. a measure space isomorphism $\theta:X\rightarrow Y$ taking $\Gamma$--orbits to $\Lambda$--orbits. For recent developments in orbit equivalence, see the surveys [Fu09],[Ga10].

In the breakthrough article [Po01], Popa studied  II$_1$ factors $M$  which admit a Cartan subalgebra satisfying both a  {\it deformation} property (in the spirit of Haagerup's property) and a {\it rigidity} property (in the spirit of the relative property (T) of Kazhdan--Margulis).  He denoted by $\Cal H\Cal T$ the class of such II$_1$ factors.  The main example of an $\Cal H\Cal T$ factor is the II$_1$ factor $M=L^{\infty}(\Bbb T^2)\rtimes$ SL$_2(\Bbb Z)$ associated with the usual action of SL$_2(\Bbb Z)$ on the 2--torus $\Bbb T^2$. 
 More generally,  if $\Gamma$ is a  group with Haagerup's property and $\Gamma\curvearrowright (X,\mu)$ is a rigid free ergodic p.m.p. action, then $M=L^{\infty}(X)\rtimes\Gamma$ is an $\Cal H\Cal T$ factor.
Recall that the action $\Gamma\curvearrowright (X,\mu)$ is  {\it rigid} if
 the inclusion $L^{\infty}(X)\subset M$ has the {\it relative property (T)}, i.e. if any sequence of unital tracial completely positive maps $\Phi_n:M\rightarrow M$ converging to the identity pointwise in $||.||_2$, must converge uniformly on the unit ball of $L^{\infty}(X)$  ([Po01]).
Here, $||.||_2$ denotes the Hilbert norm given by the trace of $M$.

The main result of [Po01] asserts that, up to unitary conjugacy, an $\Cal H\Cal T$ factor $M$ has a unique Cartan subalgebra $A$ with the relative property (T).
 The uniqueness of $A$ implies that any invariant of the inclusion $A\subset M$ is an invariant of $M$. Using this fact, Popa gave the first example of a II$_1$ factor with trivial fundamental group:  $M=L^{\infty}(\Bbb T^2)\rtimes$ SL$_2(\Bbb Z)$. Indeed, it follows that the fundamental group of $M$ is equal to the fundamental group of the orbit equivalence relation of the action SL$_2(\Bbb Z)\curvearrowright \Bbb T^2$, which is trivial by  Gaboriau's work [Ga01].

In view of [Po01] it is natural to wonder whether ${\Cal H\Cal T}$ factors have unique Cartan subalgebras. This was shown to be false in general by Ozawa and Popa in [OP08]. 
 Moreover, as noticed in [PV09] (see Section 5), their construction produces examples of  $\Cal H\Cal T$ factors that have two {\it group measure space} Cartan subalgebras.

Nevertheless, we managed to show that a large class of $\Cal H\Cal T$ factors,  which verify some rather mild assumptions (ruling out the examples from [OP08]), have a unique group measure space Cartan subalgebra.

\proclaim {Theorem 1}  Let  $\Gamma\curvearrowright (X,\mu)$ be a free ergodic rigid p.m.p. action. Assume that $\Gamma$ has positive first $\ell^2$--Betti number, $\beta_1^{(2)}(\Gamma)>0$. Denote $M=L^{\infty}(X)\rtimes\Gamma$. 

\noindent
Then $M$ has a unique group measure space Cartan subalgebra, up to unitary conjugacy. That is, if $\Lambda\curvearrowright (Y,\nu)$ is any free ergodic p.m.p. action  such that 
 $M=L^{\infty}(Y)\rtimes\Lambda$,  then we can find a unitary $u\in M$ such that 
$uL^{\infty}(X)u^*=L^{\infty}(Y)$. 
\endproclaim
Thus, if $\Gamma$ additionally has Haagerup's property, then $M$ is an $\Cal H\Cal T$ factor with a unique group measure space Cartan subalgebra. In particular, the $\Cal H\Cal T$ factor $M=L^{\infty}(\Bbb T^2)\rtimes$ SL$_2(\Bbb Z)$ has a unique group measure space decomposition. For several  concrete families of $\Cal H\Cal T$ factors with this property, see the examples below.

In their recent work [OP07], Ozawa and Popa showed that any  II$_1$ factor  $L^{\infty}(X)\rtimes\Bbb F_n$ arising from a free ergodic {\it profinite} action  of a free group $\Bbb F_n$ ($2\leqslant n\leqslant\infty$) has a unique Cartan subalgebra. Subsequently, Popa conjectured that this property should hold for {\it any} free ergodic action of $\Bbb F_n$ ([Po09]).
 Theorem 1 implies that any II$_1$ factor $L^{\infty}(X)\rtimes\Bbb F_n$ arising from a free ergodic {\it rigid} action of $\Bbb F_n$  has a unique group measure space Cartan subalgebra. 
Our result provides, thus far, the only class of actions other than [OP07] for which progress on the above conjecture has been made.

In fact, our result offers some evidence for a general conjecture which predicts that all II$_1$ factors $L^{\infty}(X)\rtimes\Gamma$ coming from free ergodic p.m.p. actions of groups $\Gamma$ with $\beta_1^{(2)}(\Gamma)>0$ have a unique Cartan subalgebra (see [Po09]).  Related to this conjecture, it has been recently shown in [CP10] (see also [Va10b]) that if $\Gamma$ additionally has a non--amenable subgroup with the relative property (T), then $L^{\infty}(X)\rtimes\Gamma$ has a unique group measure space Cartan subalgebra. 
\vskip 0.05in

We continue with several remarks on the statement of Theorem 1.

\noindent
{\it Remarks.} 
(i) We do not know whether Theorem 1 holds if instead of assuming that the action $\Gamma\curvearrowright (X,\mu)$ is rigid we only require the existence of a von Neumann subalgebra $A_0\subset L^{\infty}(X)$ such that $A_0'\cap M=L^{\infty}(X)$ and the inclusion $A_0\subset M$ has the relative property (T). When $\Gamma$ has Haagerup's property, this amounts to assuming that $A$ is an HT Cartan subalgebra rather than an HT$_s$ Cartan subalgebra ([Po01]). 
If this were the case, then [Io07, Theorem 4.3] would imply that any group $\Gamma$ with $\beta_1^{(2)}(\Gamma)>0$ admits an action whose II$_1$ factor has a unique group measure space Cartan subalgebra.

\noindent
(ii) Theorem 1 implies that the actions $\Gamma\curvearrowright (X,\mu)$ and $\Lambda\curvearrowright (Y,\nu)$ are orbit equivalent. This conclusion 
cannot be improved  to show that the groups are isomorphic and the actions are conjugate. Indeed, if $\Gamma=\Bbb F_n$, then any p.m.p. action of $\Gamma$ is orbit equivalent to actions of uncountably many non--isomorphic groups ([MS06, Theorem 2.27]).  

\noindent
(iii) Note that by [CP10, Theorem A.1] the conclusion of Theorem 1 also holds if we 
 suppose that the action $\Lambda\curvearrowright (Y,\nu)$ rather than the action $\Gamma\curvearrowright (X,\mu)$ is rigid.
\vskip 0.05in
Before providing several concrete families of actions to which Theorem 1 applies let us discuss its hypothesis. 
The study of rigid actions was initiated in [Po01] where the problem of characterizing the groups $\Gamma$ admitting a rigid action was posed. But, while this problem remains open (see [Ga08] for a partial result), several classes of rigid actions ([Po01],[Ga08],[IS10]) and an  ergodic theoretic formulation of rigidity ([Io09]) have been found.
 Recall that if $\pi:\Gamma\rightarrow\Cal O(H_{\Bbb R})$ is an orthogonal representation on a real Hilbert space $H_{\Bbb R}$, then a map $b:\Gamma\rightarrow H_{\Bbb R}$ is  a  {\it cocycle} into $\pi$ if it verifies the identity $b(gh)=b(g)+\pi(g)b(h)$, for all $g,h\in\Gamma$.
The condition $\beta_1^{(2)}(\Gamma)>0$ is equivalent to $\Gamma$ being non--amenable and having an unbounded cocycle into its left regular representation $\lambda:\Gamma\rightarrow\Cal O(\ell^2_{\Bbb R}\Gamma)$ ([BV97],[PT07]) and is satisfied by any free product group $\Gamma=\Gamma_1*\Gamma_2$ with $|\Gamma_1|\geqslant 2$ and $|\Gamma_2|\geqslant 3$. For more examples of groups with positive first $\ell^2$--Betti number, see Section 3 of [PT07]. 
\vskip 0.05in
\noindent
{\it Examples.} The following actions satisfy the hypothesis of Theorem 1:

\noindent
(i) The action  $\Gamma\curvearrowright(\Bbb T^2,\lambda^2)$, where $\Gamma<$ SL$_2(\Bbb Z)$ is a non--amenable subgroup and $\lambda^2$ is the Haar measure of $\Bbb T^2$.  

\noindent
(ii) The action $\Gamma\curvearrowright ($SL$_2(\Bbb R)/$SL$_2(\Bbb Z),m)$, where  $\Gamma$ is either a non--amenable subgroup of SL$_2(\Bbb Z)$ or a lattice of SL$_2(\Bbb R)$, and $m$ is the unique SL$_2(\Bbb R)$--invariant probability measure on SL$_2(\Bbb R)/$SL$_2(\Bbb Z).$ More generally, $\Gamma$ can be any Zariski dense countable subgroup of SL$_2(\Bbb R)$ with $\beta_1^{(2)}(\Gamma)>0$.

\noindent
(iii) Any action of the form $\Gamma\curvearrowright (G/\Lambda,m)$, where  $G$ is simple Lie group, $\Gamma<G$  is any Zariski dense countable subgroup  with $\beta_1^{(2)}(\Gamma)>0$, $\Lambda< G$ is a lattice and $m$ is the unique $G$--invariant probability measure on $G/\Lambda$. Note that by [Ku51] every semisimple Lie group $G$ contains a copy of $\Gamma=\Bbb F_2$ which is strongly dense and hence Zariski dense.

\noindent
(iv) Let $\Gamma=\Gamma_1*\Gamma_2$ be a free product group with $|\Gamma_1|\geqslant 2$ and $|\Gamma_2|\geqslant 3$. By Theorem 1.3 in [Ga08], there exists a continuum of free ergodic rigid  p.m.p. actions $\Gamma\curvearrowright (X_i,\mu_i)$, $i\in I$, such that the II$_1$ factors $L^{\infty}(X_i)\rtimes\Gamma$ are mutually non--isomorphic.
\vskip 0.02in

The groups $\Gamma$ in the examples (i)--(iv) clearly satisfy $\beta_1^{(2)}(\Gamma)>0$. The actions from (i) are rigid by [Bu91] and [Po01], while the rigidity of the actions from (ii) and (iii) is a  consequence of Theorem D in [IS10]. Note that the actions from (i)--(iii) give rise to $\Cal H\Cal T$ factors; the same is true in the case of (iv) when  $\Gamma$ has Haagerup's property.

\vskip 0.05in  The proof of Theorem 1 is based on two results that are of independent interest. The first is a structural result concerning the group measure space decompositions of II$_1$ factors $L^{\infty}(X)\rtimes\Gamma$ arising from rigid actions of groups $\Gamma$ that have an unbounded cocycle into a mixing orthogonal representation $\pi:\Gamma\rightarrow\Cal O(H_{\Bbb R})$. Recall that $\pi$ is {\it mixing} if for all $\xi,\eta\in H_{\Bbb R}$ we have that $\langle\pi(g)\xi,\eta\rangle\rightarrow 0$, as $g\rightarrow\infty$. Below we use the notation $A\prec_{M}B$ whenever ``a corner of a subalgebra $A\subset M$ can be embedded into a subalgebra $B\subset M$ inside $M$", in the sense of Popa ([Po03], see Section 1.1). This roughly means that there exists a unitary element $u\in M$ such that $uAu^*\subset B$.

\proclaim {Theorem 2} Let $\Gamma\curvearrowright (X,\mu)$ be a free ergodic rigid p.m.p. action. Assume that $\Gamma$ admits an unbounded cocycle into a mixing orthogonal representation $\pi:\Gamma\rightarrow\Cal O(H_{\Bbb R})$. Denote $M=L^{\infty}(X)\rtimes\Gamma$
and let $\Lambda\curvearrowright (Y,\nu)$ be any free ergodic p.m.p. action such that $M=L^{\infty}(Y)\rtimes\Lambda$.
For $S\subset\Lambda$, we denote by $C(S)=\{g\in\Lambda|gh=hg,\forall h\in S\}$ the centralizer of $S$ in $\Lambda$.
\vskip 0.02in
\noindent
Then we have that either

\noindent
(1) $L^{\infty}(X)\prec_{M}L^{\infty}(Y)\rtimes\Lambda_0$, for an amenable subgroup $\Lambda_0$ of $\Lambda$, or

\noindent
(2) $L^{\infty}(X)\prec_{M}L^{\infty}(Y)\rtimes(\cup_{n\geqslant 1}C(\Lambda_n))$, for a decreasing sequence $\{\Lambda_n\}_{n\geqslant 1}$  of non--amenable subgroups of $\Lambda$.
\endproclaim
The assumption that $\Gamma$ has an unbounded cocycle into a mixing representation is satisfied in particular when either $\beta_1^{(2)}(\Gamma)>0$ or $\Gamma$ has Haagerup's property.
For an outline of the proof of Theorem 2, see the beginning of Section 3. 
For now, let us mention that it uses [CP10] and, in novel fashion, ultraproduct algebras $M^{\Cal U}$ constructed from an ultrafilter $\Cal U$ over an uncountable set. 

Let us elaborate on conditions (1) and (2).
The conclusion from (1) is optimal, in the sense that it cannot be improved to deduce that $L^{\infty}(X)$ and $L^{\infty}(Y)$ are conjugate (equivalently, by [Po03], $\Lambda_0$ cannot be taken to be {\it finite}). Indeed, [OP08] provides examples of rigid actions $\Gamma\curvearrowright (X,\mu)$ of Haagerup groups $\Gamma$ whose II$_1$ factors $L^{\infty}(X)\rtimes\Gamma$ have two non--conjugate group measure space Cartan subalgebras.
Condition (2) is somewhat  imprecise in general due to our a priori lack of understanding of the subgroup structure of  $\Lambda$ and so it might seem hard to use for applications. However, in the case when $\beta_1^{(2)}(\Gamma)>0$, by using results of Chifan and Peterson [CP10] on malleable deformations  arising from cocycles into $\ell^2_{\Bbb R}\Gamma$,
we show that (2) implies (1). 

We thereby conclude that if $M=L^{\infty}(X)\rtimes\Gamma$ is as in Theorem 1 then given any  group measure space decomposition $M=L^{\infty}(Y)\rtimes\Lambda$ we can find an amenable subgroup $\Lambda_0<\Lambda$ such that $L^{\infty}(X)\prec_{M}L^{\infty}(Y)\rtimes\Lambda_0$. It follows that there is an amenable von Neumann subalgebra $N$ of $M$ such that $L^{\infty}(X)\prec_{M}N$ and $L^{\infty}(Y)\subset N$.

\vskip 0.05in
The second tool needed in the proof of Theorem 1 is a general conjugacy criterion for Cartan subalgebras which deals precisely with the last situation. 

\proclaim {Theorem 3} Let $\Gamma\curvearrowright (X,\mu)$ be a free ergodic p.m.p. action. Assume that $\beta_1^{(2)}(\Gamma)>0$ and
denote   $A=L^{\infty}(X)$, $M=A\rtimes\Gamma$. Let  $B\subset M$ be a Cartan subalgebra.

\noindent
 If there exists an amenable von Neumann subalgebra $N$ of $M$ such that $A\prec_{M}N$ and $B\subset N$,
then we can find a unitary element $u\in M$ such that $uAu^*=B$.
\endproclaim

In particular, if $A$ and $B$ generate an amenable von Neumann subalgebra of $M$, then they are unitarily conjugate.

To outline the main steps of the proof of Theorem 3 assume that $A$ and $B$ are not unitarily conjugate.  We first use the hypothesis to construct an amenable von Neumann subalgebra $P$ of $M$ such that $A\subset P$ and $B\prec_{M} P$. Secondly,  we consider the equivalence relations $\Cal R$ and $\Cal S$ on $X$ associated with the inclusions $A\subset M$ and $A\subset P$ ([FM77]). Since $B$ is regular in $M$ and has a corner which embeds into $P$ but {\it not} into $A$, we deduce that $\Cal S$ is normal in $\Cal R$, in a weak sense. Lastly, since by results of Gaboriau  an equivalence relation $\Cal R$ satisfying $\beta_1^{(2)}(\Cal R)>0$  cannot have a ``weakly normal" hyperfinite subequivalence relation ([Ga99],[Ga01]), we get a contradiction.

\vskip 0.05in
As a byproduct of the techniques developed in this paper, we also prove a rigidity result regarding the group measure space decompositions of factors $M=L^{\infty}(X)\rtimes\Gamma$ coming from actions of groups $\Gamma$ that have positive first $\ell^2$--Betti number but do not have Haagerup's property (see Theorem 6.1). We present here two interesting consequences of this result.

\proclaim {Corollary 4} Let $\Gamma$ be a countable group such that $\beta_1^{(2)}(\Gamma)\in (0,+\infty)$ and $\Gamma$ does not have Haagerup's property.
 Let $\Gamma\curvearrowright (X,\mu)$ be any free ergodic p.m.p. action.

\noindent
Then the II$_1$ factor $M=L^{\infty}(X)\rtimes\Gamma$ has trivial fundamental group, $\Cal F(M)=\{1\}$.
\endproclaim

\proclaim {Corollary 5}  Let $\Gamma$ be a countable group such that $\beta_1^{(2)}(\Gamma)>0$  and $\Gamma$ does not have Haagerup's property. Let $\Gamma\curvearrowright (X,\mu)$ be a Bernoulli action. Denote 
$M=L^{\infty}(X)\rtimes\Gamma$.

\noindent
Then $M$ has a unique group measure space Cartan subalgebra, up to unitary conjugacy.
\endproclaim

  {\it Organization of the paper.} Besides the introduction, this paper has six other sections. In Section 1, we record  Popa's intertwining technique and establish several related results. In Section 2, we review results from [CP10] on malleable deformations arising from group cocycles. Sections 3 and 4  are devoted the proofs of Theorems 2 and 3, respectively. In  Section 5 we deduce Theorem 1, while in our last section we establish Corollaries 4 and 5.

\vskip 0.05in
{\it Acknowledgment.} In the initial version of this paper, Theorems 1 and 2 were stated under the additional assumption that $\Gamma$ has Haagerup's property. I am extremely grateful to Stefaan Vaes for kindly pointing out to me that the proof of Theorem 2 can be modified to show that Theorem 2 and, consequently, Theorem 1 hold in the present generality. I would also like to thank Stefaan for allowing me to include in the text his simplified proof of Theorem 3.1.

\vskip 0.05in
{\it Added in the proof.} Very recently, Popa and Vaes proved that {\it any} II$_1$ factor arising from a free ergodic pmp action of a free group $\Gamma=\Bbb F_n$ ($2\leqslant n\leqslant\infty$) has a unique Cartan subalgebra, up to unitary conjugacy [PV11]. More generally, they showed that the same holds for any weakly amenable group $\Gamma$ with $\beta_1^{(2)}(\Gamma)>0$ [PV11] and for any hyperbolic group $\Gamma$ [PV12].

\vskip 0.2in
\head \S 1. {Preliminaries}.\endhead
\vskip 0.2in

In this paper, we work with {\it tracial von Neumann algebras} $(M,\tau)$, i.e. von Neumann algebras $M$ endowed with a faithful normal  tracial state $\tau:M\rightarrow \Bbb C$. We denote by $L^2(M)$ the completion of $M$ under the Hilbert norm $||x||_2=\tau(x^*x)^{\frac{1}{2}}$, by $\Cal U(M)$ the {\it unitary group} of $M$ and  by $(M)_1$ the {\it unit ball} of $M$, i.e. the set of $x\in M$ with $||x||\leqslant 1$. 
Given a von Neumann subalgebra $A\subset M$,  $E_A:M\rightarrow A$ denotes the {\it conditional expectation onto $A$}.

 Let us also recall the construction of the amplifications of an inclusion $A\subset M$ of a Cartan subalgebra into a II$_1$ factor. Let $t>0$.
Let $n\geqslant t$ be an integer and  $p\in D_n(\Bbb C)\otimes A$ be a projection of normalized trace $\frac{t}{n}$, where $D_n(\Bbb C)\subset\Bbb M_n(\Bbb C)$ denotes the subalgebra of diagonal matrices.
Set  $A^t:=(D_n(\Bbb C)\otimes A)p$ and  $M^t:=p(\Bbb M_n(\Bbb C)\otimes M)p$. Then the inclusion $A^t\subset M^t$, called the $t$--{\it amplification} of the inclusion $A\subset M$, is uniquely defined, up to unitary conjugacy.

\vskip 0.1in
\noindent {\bf 1.1 Popa's intertwining--by--bimodules technique}. We continue by recalling Popa's  powerful technique for conjugating subalgebras of a tracial von Neumann algebra.  Throughout this section we assume that all von Neumann algebras are separable.
 
\proclaim {Theorem 1.1 [Po03, Theorem 2.1 and Corollary 2.3]} Let $(M,\tau)$ be a tracial von Neumann algebra and  $A,N\subset M$ (possibly non--unital) von Neumann subalgebras. 
Then the following are equivalent:
\vskip 0.02in
\noindent
(1) There exist  non--zero projections $p\in A, q\in N$, a $*$--homomorphism $\psi:pAp\rightarrow qNq$ and a non--zero partial isometry $v\in qMp$ such that $\psi(x)v=vx$, for all $x\in pAp$.

\noindent
(2) There is no sequence $u_n\in\Cal U(A)$ satisfying $||E_N(au_nb)||_2\rightarrow 0$, for every $a,b\in M$.

\vskip 0.02in
\noindent
If these  equivalent conditions hold true,  we say that {\it a corner of $A$ embeds into $N$ inside $M$} and write $A\prec_{M}N$.
\endproclaim

\vskip 0.05in
\noindent
{\it Remark 1.2}. Assume that $N_1,..,N_k\subset M$ are von Neumann subalgebras such that $A\nprec_{M}N_i$, for all $i\in\{1,..,k\}$. Then  we can find a sequence $u_n\in\Cal U(A)$ such that $||E_{N_i}(au_nb)||_2\rightarrow 0$, for all $a,b\in M$ and every $i\in\{1,..,k\}$. 

To see this, identify $A$ with the diagonal subalgebra $\{(x\oplus..\oplus x)|x\in A\}$ of $\tilde M=\bigoplus_{i=1}^kM$ and let $N=\bigoplus_{i=1}^kN_i\subset\tilde M$.   Since $A\nprec_{M}N_i$, for all $i$,  the first part of Theorem 1.1 implies that $A\nprec_{\tilde M}N$. Thus, by part (2) of Theorem 1.1 we can find $u_n\in\Cal U(A)$ such that $||E_{N}(au_nb)||_2\rightarrow 0$, for all $a,b\in\tilde M$. This sequence clearly satisfies our claim.

\vskip 0.05in

Next, we record several useful related results. The first, due to Popa, asserts that for Cartan subalgebras, ``embedability of a corner" is equivalent to unitary conjugacy.

\proclaim {Lemma 1.3 [Po01, Theorem A.1.]} Let $M$ be a II$_1$ factor and $A,B\subset M$  two Cartan subalgebras. If $A\prec_{M}B$, then we can find  $u\in \Cal U(M)$ such that $uAu^*=B$.\endproclaim

\proclaim {Lemma 1.4 [PP86]} 
Let $(M,\tau)$ be a tracial von Neumann algebra and $A,N\subset M$ two von Neumann subalgebras. If $A\nprec_{M}N$, then for every $\varepsilon>0$ we can find a projection $e\in A$ such that $||E_N(e)||_2<\varepsilon ||e||_2$.
\endproclaim
\noindent
{\it Proof.} It is easy to see that $A$ and $N$ can be assumed unital. Let $\langle M,e_N\rangle$ be {\it Jones' basic construction} of the inclusion $N\subset M$ endowed with its natural semi--finite trace $Tr$.  If $A\nprec_{M}N$,  by Theorem 2.1 in [Po03], $A'\cap \langle M,e_N\rangle$ contains no projections of finite trace. Let $\varepsilon>0$. By applying Lemma 2.3. in [PP86], we can find projections $e_1,..,e_n\in M$ such that $\sum_{i=1}^ne_i=1$ and $||\sum_{i=1}^ne_ie_Ne_i||_{2,Tr}<\varepsilon$. Since $||\sum_{i=1}^{n}e_ie_Ne_i||_{2,Tr}^2=\sum_{i=1}^n||E_N(e_i)||_2^2$, we can find $i$ such that $e=e_i$ satisfies the conclusion.\hfill$\blacksquare$

\proclaim {Lemma 1.5}  Let $(M,\tau)$ be a tracial von Neumann algebra and $A,N\subset M$ two von Neumann subalgebras. Assume that $A$ is maximal abelian in $M$ and $A\prec_{M}N$.  

\vskip 0.03in
\noindent
Then there exist projections $p\in A, q\in N$, a $*$--homomorphism $\psi:Ap\rightarrow qNq$ and a non--zero partial isometry $v\in qMp$ such that $\psi(x)v=vx$, for all $x\in Ap$, and $\psi(Ap)$ is maximal abelian in $qNq$.
\endproclaim

\noindent
{\it Proof}. By the hypothesis we can find projections $p\in A,q\in N$, a $*$--homomorphism $\psi:Ap\rightarrow qNq$ and a non--zero partial isometry $v\in qMp$ such that $\psi(x)v=vx$, for all $x\in Ap$, $v^*v=p$ and $q':=vv^*\in\psi(Ap)'\cap qMq$. 
After replacing $q$ with a subprojection, we may assume that $q$ is the support projection of $E_N(q')$ and that $cq\leqslant E_N(q')\leqslant Cq$, for some $c,C>0$. 
Denote $\Cal A=\psi(Ap)'\cap qNq$.

\vskip 0.05in
\noindent
{\bf Claim.}  $\psi(Ap)q_0$ is maximal abelian in $q_0Nq_0$, for some non--zero projection $q_0\in\Cal A$.

\vskip 0.05in

Assuming the claim, define $\psi_0:Ap\rightarrow q_0Nq_0$  by $\psi_0(x)=\psi(x)q_0$ and let $v_0=q_0v$. Since $\psi_0(x)v_0=v_0x$ for all $x\in Ap$ the claim implies the lemma.

\vskip 0.05in

Now, the claim follows from Step 2 in the proof of [Po01, Theorem A.2.]. For completeness, we provide a proof.
\vskip 0.05in
\noindent
{\it Proof of the claim}.
 Since $\psi(Ap)q'=vApv^*$ and $A$ is maximal abelian in $M$, we get that $q'(\psi(Ap)'\cap qMq)q'=\psi(Ap)q'$. Fix a projection $e\in\Cal A$ and let $f\in\psi(Ap)$, $0\leqslant f\leqslant q$, such that $q'eq'=fq'$. Since $fq=f\in \psi(Ap)\subset N$ and $E_N(q')\geqslant cq$, we have that $||e||_2\geqslant ||fq'||_2=\tau(f^2q')^{\frac{1}{2}}=\tau(f^2E_N(q'))^{\frac{1}{2}}\geqslant c^{\frac{1}{2}}\tau(f^2)^{\frac{1}{2}}=c^{\frac{1}{2}}||f||_2$.

 Further, since $e,f\in N$ and $f\in\psi(Ap)$, we have that

$$ ||eq'e||_2^2=\tau(efq')=\tau(efE_N(q'))\leqslant C\tau(ef)\leqslant C||E_{\psi(Ap)}(e)||_2||f||_2\leqslant\tag 1.a$$ $$ Cc^{-\frac{1}{2}}||E_{\psi(Ap)}(e)||_2||e||_2.$$

On the other hand, since $e\in N$ and $E_N(q')\geqslant cq$, we get that

$$||eq'e||_2\geqslant ||E_N(eq'e)||_2=||eE_N(q')e||_2\geqslant c||e||_2\tag 1.b$$

Combining (1.a) and (1.b) yields that 
$||E_{\psi(Ap)}(e)||_2\geqslant C^{-1}c^{\frac{3}{2}}||e||_2,$ for any projection $e\in\Cal A.$
Since $\psi(Ap)$ is abelian, Lemma 1.4 and Theorem 1.1 imply that $\Cal A$ is of type $\text{I}_{fin}$. Hence, if we denote by $\Cal Z$ the center of $\Cal A$, then we can find a non--zero projection $q_1\in\Cal A$ such that $q_1\Cal Aq_1=\Cal Zq_1$. The last inequality and Lemma 1.4  also imply that $\Cal Zq_1\prec_{\Cal A}\psi(Ap)$. Thus, $\psi(Ap)q_0=\Cal Zq_0=q_0\Cal Aq_0$, for non--zero projection $q_0\in\Cal Zq_1$. This finishes the proof of the claim and of the lemma.
\hfill$\blacksquare$

\vskip 0.02in

\proclaim {Lemma 1.6} Let $(M,\tau)$ be a tracial von Neumann algebra, $N\subset M$ a von Neumann subalgebra and $q\in M$  a projection.  Let $q_0$ be the support projection of $E_N(q)$.

\vskip 0.02in
\noindent (1) If we denote by $P\subset q_0Nq_0$ the von Neumann algebra generated by $E_N(qMq)$, then $pNp\prec_{N}Pp$, for every non--zero projection $p\in P'\cap q_0Nq_0$.
\vskip 0.02in
\noindent (2)  If we denote by $Q\subset qMq$ the von Neumann algebra generated by $qNq$, then $pNp\prec_{M}Q$, for every non--zero projection $p\in q_0Nq_0$.
\endproclaim

\noindent
{\it Proof.}  Using functional calculus for the positive operator $E_N(q)$, we define $q_t=1_{[t,1]}(E_N(q))$, for every $t\in [0,1]$. Then $q_t\in P$ and $||q_t-q_0||_2\rightarrow 0$, as $t\rightarrow 0$.
\vskip 0.05in
\noindent
(1) Let $p\in P'\cap q_0Nq_0$. Then $p_t=pq_t$ is a projection and $||p_t-p||_2\rightarrow 0$, as $t\rightarrow 0$. In order to get the conclusion, it suffices to prove that $p_tNp_t\prec_{N}Pp$, for all $t>0$.
 Let $e\in p_tNp_t$ be a projection. Since $e=ep\in N$ and $pE_N(qeq)\in Pp$ we have that 

$$||eqe||_2^2=\tau(epqeq)=\tau(epE_N(qeq))=\tau(E_{Pp}(e)pE_N(qeq))\leqslant ||E_{Pp}(e)||_2||e||_2\tag 1.c$$

On the other hand, since $e=p_te$ and $E_N(q)p_t\geqslant tp_t\geqslant 0,$ we get
$$||eqe||_2^2\geqslant ||E_N(eqe)||_2^2=||eE_N(q)e||_2^2=||eE_N(q)p_te||_2^2\geqslant t^2||e||_2^2\tag 1.d $$

Combining (1.c) and (1.d) yields that $||E_{Pp}(e)||_2\geqslant t^2||e||_2$, for all projections $e\in p_tNp_t$. Then Lemma 1.4 implies that $p_tNp_t\prec_{N}Pp$, as claimed.

\vskip 0.05in
\noindent (2). Since $||q_t-q_0||_2\rightarrow 0$, we may assume that $p\leqslant q_t$, for some $t>0$. 
 Let $e\in pNp$ be a projection. Then $qeq\in Q$, hence $\tau(eqeq)=\tau(E_{Q}(e)qeq)\leqslant ||E_{Q}(e)||_2||e||_2$. On the other hand, since $E_N(eqe)=eE_N(q)e=eE_N(q)q_te\geqslant te$, as in (1.d) we get that $\tau(eqeq)=||eqe||_2^2\geqslant t^2||e||_2^2$.

The last two inequalities together imply that $||E_{Q}(e)||_2\geqslant t^2||e||_2$, for any projection $e\in pNp$. By applying Lemma 1.4 we obtain that $pNp\prec_{M}Q$.
\hfill$\blacksquare$

\vskip 0.1in
\noindent {\bf 1.2 Equivalence relations from Cartan subalgebras.} Consider a standard probability space $(X,\mu)$. 
 A Borel equivalence relation $\Cal R\subset X^2$  is called {\it countable}, {\it measure preserving} if it is induced by a measure preserving action of a countable group on $(X,\mu)$ ([FM77]). 
We denote by $[\Cal R]$ (the {\it full group} of $\Cal R$) the group of Borel automorphisms  $\theta$ of $X$ such that $\theta(x)\Cal Rx$, for almost all $x\in X$. Also, we denote by $[[\Cal R]]$ (the {\it full pseudogroup} of $\Cal R$) the set of Borel isomorphisms $\theta:Y\rightarrow Z$ satisfying $\theta(x)\Cal Rx$, for almost all $x\in Y$, where $Y,Z\subset X$ are Borel sets.

Next, we recall the construction of equivalence relations coming from Cartan subalgebra inclusions.
 Let $(M,\tau)$ be a separable tracial von Neumann algebra with a Cartan subalgebra $A$. Identify $A$ with $L^{\infty}(X)$, where $(X,\mu)$ is a standard probability space. Every $u\in \Cal N_{M}(A)$ defines an automorphism $\theta_u$ of $(X,\mu)$ by $a\circ\theta_u=u^*au$, for $a\in A$. Let $\Gamma<\Cal N_M(A)$ be a countable, $||.||_2$--dense subgroup. The {\it equivalence relation of the inclusion} $(A\subset M)$, denoted $\Cal R_{(A\subset M)}$, is given by $x\sim y$ iff $x=\theta_u(y)$, for some $u\in\Gamma$.

Note that $\Cal R_{(A\subset M)}$ is countable, measure preserving and does not depend on the choice of $\Gamma$. The latter is a consequence of the following fact: if $u\in\Cal N_{M}(A)$ and $u_n\in\Gamma$ are such that $||u_n-u||_2\rightarrow 0$, then $\mu(\{\theta_{u_n}=\theta_u\})\rightarrow 0$ and thus $\theta_u\in [\Cal R_{(A\subset M)}]$.

For later reference, we fix the following notation. If $\theta:Y\rightarrow Z$ belongs to $[[\Cal R_{(A\subset M)}]]$, then we can find a partial isometry $u_{\theta}\in M$ which ``implements" $\theta$: $u_{\theta}u_{\theta}^*=1_{Z}$, $u_{\theta}^*u_{\theta}=1_Y$ and $u_{\theta}^*au_{\theta}=(a\circ\theta)1_Y$, for all $a\in A$.

\vskip 0.05in

The next lemma is the analogue of Popa's intertwining technique (Theorem 1.1) for equivalence relations. Note that it generalizes part of Theorem 2.5. in [IKT08].

\proclaim {Lemma 1.7} Let $\Cal R$ be a countable, measure preserving equivalence relation on a probability space $(X,\mu)$. Let $\Cal S,\Cal T$ be two subequivalence relations. 

\noindent
Define $\varphi_{\Cal S}:[\Cal R]\rightarrow [0,1]$ by $\varphi_{\Cal S}(\theta)=\mu(\{x\in X|\theta(x)\Cal Sx\})$. 
 Assume that there is no sequence $\{\theta_n\}_{n\geqslant 1}\subset [\Cal T]$ such that $\varphi_{\Cal S}(\psi\theta_n\psi')\rightarrow 0$, for all $\psi,\psi'\in [\Cal R]$.

\vskip 0.03in
\noindent
Then we can find $\theta\in[[\Cal R]]$, with $\theta:Y\rightarrow Z$, and $k\geqslant 1$ such that every $(\theta\times\theta)(\Cal T_{|Y})$--class is contained in the union of at most $k$ $\Cal S_{|Z}$--classes.

\endproclaim

\noindent
{\it Proof.} We first claim that there are $\psi_1,..,\psi_k,\psi_1',..,\psi_k'\in [\Cal R]$ and $c>0$ such that $$\sum_{i,j=1}^{k}\varphi_{\Cal S}(\psi_i\theta\psi_j')\geqslant c,\hskip 0.05in\forall \theta\in [\Cal T]\tag 1.e$$

Assume by contradiction that this is false. Fix two sequences $\{\psi_i\}_{i\geqslant 1},\{\psi_j'\}_{j\geqslant 1}\subset [\Cal R]$
which are dense with respect to the metric $d(\theta_1,\theta_2)=\mu(\{\theta_1\not=\theta_2\})$. Then by our assumption, we can find a sequence $\{\theta_n\}_{n\geqslant 1}\subset [\Cal T]$ such that $\varphi_{\Cal S}(\psi_i\theta_n\psi_j')\rightarrow 0$, for all $i,j\geqslant 1$. Using the density of $\{\psi_i\}_{i\geqslant 1}$ and $\{\psi_j'\}_{j\geqslant 1}$, it follows that $\varphi_{\Cal S}(\psi\theta_n\psi')\rightarrow 0$, for all $\psi,\psi'\in [\Cal R]$, contradicting the hypothesis.

In the rest of the proof we follow closely Section 2 of [IKT08].
 First, we may assume that every $\Cal R$--class contains infinitely many $\Cal S$--classes. Thus, we can find a sequence of Borel functions $C_n:X\rightarrow X$ such that $C_0=$ id and for a.e. $x\in X$, $\{C_n(x)\}_{n\geqslant 0}$ is a transversal for the $\Cal S$--classes contained in the $\Cal R$--class of $x$. 

Denote by $S(\Bbb N)$ be the symmetric group of $\Bbb N$ and by $\rho$ the counting measure on $\Bbb N$.
As in Section 2 of [IKT08], define the cocycle $w:\Cal R\rightarrow S(\Bbb N)$  by $w(x,y)(m)=n\Longleftrightarrow (C_m(x),C_n(y))\in\Cal S.$
Further, define the group morphism $\pi:[\Cal R]\rightarrow$ Aut$(X\times\Bbb N,\mu\times\rho)$ by the formula $\pi(\theta)(x,m)=(\theta(x),w(\theta(x),x)(m)),$ for all $\theta\in [\Cal R]$ and $(x,m)\in X\times\Bbb N$. Denote also by $\pi$ the associated unitary representation of $[\Cal R]$ on $\Cal H=L^2(X\times \Bbb N)$.

Set $\xi_0=1_{X\times\{0\}}\in\Cal H$. Then $\varphi_{\Cal S}(\theta)=\langle\pi(\theta)(\xi_0),\xi_0\rangle$, for all $\theta\in [\Cal R]$.  Thus (1.e) rewrites as $\sum_{i,j=1}^{k}\langle\pi(\theta)(\pi(\psi_j')(\xi_0)),\pi(\psi_i^{-1})(\xi_0)\rangle\geqslant c$, for all $\theta\in [\Cal T]$. This implies that the restriction of $\pi$ to $[\Cal T]$ is not weakly mixing.
Let  $\xi\in \Cal H\overline{\otimes}\Cal H\cong L^2((X\times\Bbb N,\mu\times\rho)^2)$ be a non--zero $(\pi\otimes\pi)([\Cal T])$--invariant vector. 

\vskip 0.05in
\noindent
{\bf Claim}. We have that $(\pi(\theta)\otimes 1)(\xi)=\xi$, for all $\theta\in [\Cal T]$.

\vskip 0.05in
\noindent
{\it Proof of the claim.}
Let $\theta\in [\Cal T]$. Then we can find a sequence $\theta_n\in [\Cal T]$ such that for almost every $(x,y)\in X^2$ we may find $n\geqslant 1$ satisfying $\theta(x)=\theta_n(x)$ and $y=\theta_n(y)$. Since $(\pi(\theta_n)\otimes\pi(\theta_n))(\xi)=\xi$ it follows easily that
  $(\pi(\theta)\otimes 1)(\xi)=\xi$.

To construct a sequence as above, let $n\geqslant 1$ and consider a partition $A_1,..,A_n$ of $X$ with $\mu(A_i)=\frac{1}{n}$. For $1\leqslant i\leqslant n$, let $\theta_{i,n}\in [\Cal T]$ such that $\theta_{i,n}(x)=\theta(x)$, for $x\in A_{i,n}$ and $\theta_{i,n}(y)=y$, for $y\in X\setminus (A_{i,n}\cup\theta(A_{i,n}))$. Let $Y_n$ be the set of $(x,y)\in X^2$ for which we may find $i\in\{1,..,n\}$ with $\theta(x)=\theta_{i,n}(x)$ and $y=\theta_{i,n}(y)$. Since $Y_n$ contains 
$A_{i,n}\times (X\setminus (A_{i,n}\cup\theta(A_{i,n})),$ for all $i$, we get that $(\mu\times\mu)(Y_n)\geqslant 1-\frac{2}{n}$. Thus $\cup_{n\geqslant 1}Y_n=X^2$, implying that the sequence $\{\theta_{i,n}\}_{1\leqslant i\leqslant n<\infty}$ verifies the desired conditions.
\hfill$\square$
\vskip 0.02in

The claim implies that we can find a non--zero $\pi([\Cal T])$--invariant vector $\eta\in\Cal H$.
For $x\in X$, let $N_x=\{n\in\Bbb N|\hskip 0.02in|\eta(x,n)|$ is maximal among all $|\eta(x,i)|$, $i\in\Bbb N\}$. 
Since $\eta$ is $\pi([\Cal T])$--invariant it follows that $w(y,x)N_x=N_y$, for almost all $(x,y)\in\Cal T$.
Since $\eta\in L^2(X\times\Bbb N)$, we can find $\kappa\geqslant 1$ and a set $X_0\subset X$ of positive measure such that $|N_x|=\kappa$, for every $x\in X_0$. Enumerate $N_x=\{n_{1,x},..,n_{\kappa,x}\}$ and let $n_x=n_{1,x}$. 

Define the equivalence relation $\Cal T_0$ on $X_0$ as the set of $(x,y)\in\Cal T\cap (X_0\times X_0)$ such that $w(y,x)n_{i,x}=n_{i,y}$, for all $1\leqslant i\leqslant \kappa$. Since for  all $(x,y)\in\Cal T$ we can find a permutation $\pi$ of $\{1,..,\kappa\}$ such that $n_{i,y}=w(y,x)n_{\pi(i),x}$, it follows that every $\Cal T_{|X_0}$--class contains at most $k:=\kappa!$ $\Cal T_0$--classes.

Now, for almost all $(x,y)\in\Cal T_0$ we have $w(y,x)n_{x}=n_{y}$, thus $(C_{n_{x}}(x),C_{n_{y}}(y))\in\Cal S$. 
Let $Y\subset X_0$ be a set of positive measure such that the map $Y\ni x\rightarrow \theta(x)=C_{n_x}(x)$ is 1--1. It follows that $\theta:Y\rightarrow Z=\theta(Y)$ belongs to $[[\Cal R]]$ and  $(\theta\times\theta)({\Cal T_0}_{|Y})\subset \Cal S_{|Z}$. Since every $\Cal T_{|Y}$--class is contained in the union of at most $k$ ${\Cal T_0}_{|Y}$--classes, we are done.
\hfill$\blacksquare$

\vskip 0.05in

\proclaim {Lemma 1.8} Let $(M,\tau)$ be a separable tracial von Neumann algebra, $A\subset M$ a Cartan subalgebra and $N,P\subset M$ von Neumann subalgebras containing $A$. Identify $A=L^{\infty}(X),$ where $(X,\mu)$ is a probability space. Let $\Cal R=\Cal R_{(A\subset M)}$, $\Cal S=\Cal R_{(A\subset N)}$ and $\Cal T=\Cal R_{(A\subset P)}$. 
\vskip 0.03in
\noindent
Then  $P\prec_{M}N$ if and only if  we can find $\theta\in[[\Cal R]]$, with $\theta:Y\rightarrow Z$, and $k\geqslant 1$ such that every $(\theta\times\theta)(\Cal T_{|Y})$--class is contained in the union of at most $k$ $\Cal S_{|Z}$--classes. 
\endproclaim
\noindent {\it Proof.} The ``if" part follows easily and we leave its proof to the reader. For the ``only if" part assume that we cannot find $\theta\in [[\Cal R]]$ and $k\geqslant 1$ as above. Lemma 1.7 then provides a sequence $\theta_n\in [\Cal T]$ such that $\varphi_{\Cal S}(\psi\theta_n\psi')\rightarrow 0$, for all $\psi,\psi'\in [\Cal R]$. We claim that $||E_N(xu_{\theta_n}y)||_2\rightarrow 0$, for all $x,y\in M$. Since $u_{\theta_n}\in\Cal U(P)$, it follows that $P\nprec_{M}N$. Thus, the claim finishes the proof of the ``only if" part.

Since $E_P$ is $A$--bimodular, by Kaplansky's theorem it suffices to prove the claim for $x=u_{\psi}$ and $y=u_{\psi'}$, where $\psi,\psi'\in [\Cal R]$. In this case, $||E_N(u_{\psi}u_{\theta_n}u_{\psi'})||_2=\sqrt{\varphi_{\Cal S}(\psi\theta_n\psi')}\rightarrow 0$, as claimed.  \hfill$\blacksquare$

\vskip 0.2in
\head \S 2. {Deformations from group cocycles.}\endhead

\vskip 0.2in 
Let  $(A,\tau)$ be a tracial von Neumann algebra, $\Gamma\curvearrowright A$ be a trace preserving action and set $M=A\rtimes\Gamma$. 
Let $\pi:\Gamma\rightarrow\Cal O(H_{\Bbb R})$ be an orthogonal representation, where $H_{\Bbb R}$ is a separable real Hilbert space. 
Given a cocycle $b:\Gamma\rightarrow H_{\Bbb R}$, Sinclair constructed a {\it malleable deformation} in the sense of Popa, i.e. a tracial von Neumann algebra $\tilde M\supset M$ and a 1--parameter group of automorphisms $\{\alpha_t\}_{t\in\Bbb R}$ of $\tilde M$ such that $||\alpha_t(x)-x||_2\rightarrow 0$ for all $x\in\tilde M$  (see [Si10, Section 3] and [Va10b, Section 3.1]).

To recall this construction, fix an orthonormal basis $\Cal B\subset H_{\Bbb R}$ and let $(X,\mu)=\prod_{v\in \Cal B}(\Bbb R,\mu_0)_{v}$, where $d\mu_0=\frac{1}{\sqrt{2\pi}}\exp(-\frac{x^2}{2})dx$ is the Gaussian measure on $\Bbb R$. 

Next,  for every $\xi=\sum_{v\in\Cal B}c_{v}v\in H_{\Bbb R}$ (with $c_v\in\Bbb R$) we define a unitary $\omega(\xi)\in L^{\infty}(X)$ by letting $\omega(\xi)(x)=\exp(\sqrt{2}i\sum_{v\in\Cal B}c_{v}x_{v})$, for all $x=(x_v)_{v}\in X.$
Then  $\omega(\xi+\eta)=\omega(\xi)\omega(\eta)$, $\omega(\xi)^*=\omega(-\xi)$ and $\tau(\omega(\xi))=\exp(-||\xi||^2)$, for all $\xi,\eta\in H_{\Bbb R}$.

Define $D\subset L^{\infty}(X)$ to be the von Neumann algebra generated by $\{\omega(\xi)|\xi\in H_{\Bbb R}\}$ and let $\tau$ be the trace on $D$ given by integration against $\mu$. Consider the Gaussian action $\Gamma\curvearrowright^{\sigma} D$ which on the generating functions $\omega(\xi)$ is given by $\sigma_g(\omega(\xi))=\omega(\pi(g)(\xi))$. Finally, let 
$\Gamma\curvearrowright D\overline{\otimes}A$ be the diagonal action and define $\tilde M=(D\overline{\otimes}A)\rtimes\Gamma$.

It follows that the formula   $$\alpha_t(u_g)=(\omega(tb(g))\otimes 1)u_g \hskip 0.05in\text{for all}\hskip 0.05in g\in\Gamma\hskip 0.05in\text{and}\hskip 0.05in \alpha_t(x)=x\hskip 0.05in\text{for all}\hskip 0.05in x\in D\overline{\otimes}A$$ gives a 1--parameter group of automorphisms $\{\alpha_t\}_{t\in\Bbb R}$ of $\tilde M$.
Note that $\alpha_t\rightarrow id$ in the pointwise $||.||_2$--topology: $||\alpha_t(x)-x||_2\rightarrow 0$, for all $x\in\tilde M$. 
Given $S\subset \tilde M$ we say that $\alpha_t\rightarrow id$ {\it uniformly on $S$} if $\sup_{x\in S}||\alpha_t(x)-x||_2\rightarrow 0$, as $t\rightarrow 0$. 
\vskip 0.05in

Next, we recall several results concerning the deformations $\{\alpha_t\}_{t\in\Bbb R}$ that we will subsequently need.

\proclaim {Lemma 2.1}  
 If $\alpha_t\rightarrow id$ uniformly on $(pMp)_1$, for some non--zero projection $p\in M$, then $b$ is a bounded cocycle.
\endproclaim
\noindent
{\it Proof.}  If $\alpha_t\rightarrow id$ uniformly on $(pMp)_1$, then $\alpha_t\rightarrow id$ uniformly on $(Mz)_1$, where $z$ is the central support of $p$ in $M$. Therefore $\tau(\alpha_t(u_g)u_g^*z)\rightarrow \tau(z)$, uniformly in $g\in\Gamma$.
Since $E_{M}(\alpha_t(u_g))=\exp(-t^2||b(g)||^2)u_g$, we deduce that $\exp(-t^2||b(g)||^2)\rightarrow 1$, uniformly in $g\in\Gamma$. This implies that $b$ is bounded.\hfill$\blacksquare$

\proclaim {Lemma 2.2 [Po06b]} Let $p\in M$ be a projection and $B\subset pMp$ be a von Neumann algebra. 
If $\pi$ is weakly contained in the left regular representation of $\Gamma$ and $B$ has no amenable direct summand, then $\alpha_t\rightarrow id$ uniformly on $(B'\cap pMp)_1$.
\endproclaim
\noindent {\it Proof.} This is a direct consequence of Popa's spectral gap argument. For the reader's convenience let us sketch a proof. Since $\pi$ is weakly contained in the left regular representation of $\Gamma$, the $M$--$M$ bimodule $L^2(\tilde M)\ominus L^2(M)$ is weakly contained in the $M$--$M$ bimodule $(L^2(M)\overline{\otimes}L^2(M))^{\oplus\infty}$ (see e.g. [Va10b, Lemma 3.5]). 

Fix $\varepsilon>0$. Since $B$ has no amenable direct summand,  the proof of [Po06b, Lemma 2.2] shows that   we can find $b_1,..,b_n\in B$ and $\delta>0$ such that if $x\in p\tilde Mp$ satisfies $||x||\leqslant 1$ and $||[x,b_i]||_2\leqslant\delta$, for all $i\in\{1,..,n\}$, then $||x-E_M(x)||_2\leqslant\varepsilon.$

Next, we use Popa's spectral gap argument (see the proof of [Po06b, Theorem 1.1]). Choose $t_0$ such that for all $|t|\leqslant t_0$ we have that $||\alpha_{-t}(b_i)-b_i||_2\leqslant\frac{\delta}{4}$ and $||\alpha_{-t}(p)-p||_2\leqslant\min{\{\frac{\delta}{8},\varepsilon\}}$. Fix $x\in (B'\cap pMp)_1$ and $t$ with $|t|\leqslant t_0$. Since $[b_i,pxp]=0$, we get that $$||[b_i,p\alpha_t(x)p]||_2=||[\alpha_{-t}(b_i),\alpha_{-t}(p)x\alpha_{-t}(p)]||_2\leqslant$$ $$ 2||\alpha_{-t}(b_i)-b_i||_2+4||\alpha_{-t}(p)-p||_2\leqslant\delta, \hskip 0.05in\forall i\in\{1,..,n\}.$$

It follows that $||p\alpha_t(x)p-E_{M}(p\alpha_t(x)p)||_2\leqslant\varepsilon$. Since $||\alpha_t(x)-p\alpha_t(x)p||_2\leqslant 2||\alpha_t(p)-p||_2\leqslant 2\varepsilon$, we get that $||\alpha_t(x)-E(\alpha_t(x))||_2\leqslant 3\varepsilon$. Finally, [Va10b, Lemma 3.1] implies that $||\alpha_t(x)-x||_2\leqslant 3\sqrt{2}\varepsilon$. Since this happens for all $t\in\Bbb R$ with $|t|\leqslant t_0$ and every $x\in (B'\cap pMp)_1$, we are done.\hfill$\blacksquare$

\vskip 0.05in
Let $B\subset M$ be a von Neumann subalgebra. 
Peterson [Pe06, Theorem 4.5] and Chifan and Peterson  [CP10, Theorem 2.5]  proved that if $\alpha_t\rightarrow id$ uniformly on $(B)_1$ and $B\nprec_{M}A$ then $\alpha_t\rightarrow id$ uniformly on  $\Cal N_{M}(B)$.

\proclaim {Theorem 2.3 [Pe06] and [CP10]} Assume that $\pi$ is mixing. Let $p\in M$ be a projection and $B\subset pMp$ be a von Neumann subalgebra. Denote by $P$ the von Neumann algebra generated by the normalizer of $B$ inside $pMp$. 

\noindent
 If $\alpha_t\rightarrow id$ uniformly on $(B)_1$ and $B\nprec_{M}A$, then $\alpha_t\rightarrow id$ uniformly on $(P)_1$. 

\endproclaim

 Conversely, Chifan and Peterson proved in [CP10, Theorem 3.2] that if $B$ is abelian and $\alpha_t\rightarrow id$ uniformly on  a sequence $\{u_k\}_{k\geqslant 1}\subset\Cal N_{M}(B)$  which ``converges weakly to 0 relative to $A$", then $\alpha_t\rightarrow id$ on $(B)_1$. More generally, we have

\proclaim {Theorem 2.4 [CP10]} Assume that $\pi$ is mixing. Let $p\in M$ be a projection and $B\subset pMp$ be an abelian von Neumann subalgebra. 
Assume that  we can find a  net $(u_j)_{j\in J}$ of unitary elements in $pMp$ which normalize $B$  such that  

$\bullet$ $\alpha_t\rightarrow id$ uniformly on the tail of $(u_j)_{j\in J}$ and 

$\bullet$ $\lim_j||E_A(xu_jy)||_2=0$, for all $x,y\in M$.
\vskip 0.02in
\noindent Then $\alpha_t\rightarrow id$ uniformly on $(B)_1$.
\endproclaim

Here, following [Va10b], we say that $\alpha_t\rightarrow id$ uniformly on the tail of $(u_j)_{j\in J}$ if for all $\varepsilon>0$ we can find $j_0\in J$ and $t_0>0$ such that $||\alpha_t(u_j)-u_j||_2\leqslant\varepsilon$, for all $j\geqslant j_0$ and every $|t|\leqslant t_0$.

Theorems 2.3 and 2.4 were proved in [Pe06] and [CP10] using Peterson's technique of unbounded derivations [Pe06]. For proofs using the 1--parameter group of automorphisms $\{\alpha_t\}_{t\in\Bbb R}$, see Vaes's paper [Va10b, Theorems 3.9 and 4.1].

\vskip 0.05in

We end this section with two  facts about cocycles  (see e.g. [Pe06, Section 4]), which can be viewed as group--theoretic counterparts of 2.2 and 2.3:

\proclaim {Lemma 2.5} Let $\pi:\Gamma\rightarrow\Cal O(H_{\Bbb R})$ be an orthogonal representation and $b:\Gamma\rightarrow H_{\Bbb R}$ be a cocycle for $\pi$. Let $\Gamma_0<\Gamma$ be a subgroup.

\vskip 0.02in
\noindent
(1) If $\pi$ is weakly contained in the left regular representation of $\Gamma$ and $\Gamma_0$ is non--amenable, then the restriction of $b$ to the centralizer of $\Gamma_0$ is bounded.

\noindent
(2) Assume that $\pi$ is mixing and that $b(g)=\lambda(g)\xi-\xi$, for all $g\in\Gamma_0$, for some $\xi\in\ell^2\Gamma$. Let $h\in\Gamma$. If  $h\Gamma_0h^{-1}\cap \Gamma_0$ is infinite, then $b(h)=\lambda(h)\xi-\xi$. 
\endproclaim
 
\noindent
{\it Proof}. (1) Since $\Gamma_0$ is non--amenable, the restriction of $\pi$ to $\Gamma_0$ does not have almost invariant vectors. Hence we can find $g_1,..,g_n\in\Gamma_0$ such that $||\xi||\leqslant\sum_{i=1}^n||\pi(g_i)\xi-\xi||$, for all $\xi\in\ell^2\Gamma$. It follows that if $g\in\Gamma$ is in the centralizer of $\Gamma_0$, then $||b(g)||\leqslant\sum_{i=1}^n||\pi(g_i)b(g)-b(g)||=\sum_{i=1}^n||\pi(g)b(g_i)-b(g_i)||\leqslant 2\sum_{i=1}^n||b(g_i)||$.

\vskip 0.02in
\noindent
 (2) Define a new cocycle $\tilde b$ by letting $\tilde b(g)=b(g)-(\pi(g)\xi-\xi)$, for $g\in\Gamma$. Then $\tilde b(g)=0$, for all $g\in\Gamma_0$. Let $h\in\Gamma$ with $h\Gamma_0h^{-1}\cap \Gamma_0$ infinite and fix $g\in h\Gamma_0h^{-1}\cap \Gamma_0$. Let $k\in\Gamma_0$ such that $gh=hk$. Since $\tilde b(g)=\tilde b(k)=0$, we get that $\pi(g)\tilde b(h)=\tilde b(h)$, for all $g\in h\Gamma_0h^{-1}\cap \Gamma_0$. Since $\pi$ is a mixing representation it follows that $\tilde b(h)=0$.\hfill$\blacksquare$
\vskip 0.2in

\vskip 0.2in 
\head \S 3. {A structural result for group measure space decompositions}.\endhead

\vskip 0.2in

In this section we prove the following generalization of Theorem 2:

\proclaim {Theorem 3.1} Let $\Gamma\curvearrowright (X,\mu)$ be a free ergodic p.m.p. action and denote $A=L^{\infty}(X)$ and  $M=A\rtimes\Gamma$. Assume that $\Gamma$ admits an unbounded cocycle $b:\Gamma\rightarrow H_{\Bbb R}$ into a mixing orthogonal representation $\pi:\Gamma\rightarrow\Cal O(H_{\Bbb R})$.

\noindent 
Assume that $M^t=L^{\infty}(Y)\rtimes\Lambda$, for a free ergodic p.m.p. action $\Lambda\curvearrowright (Y,\nu)$ and $t>0$. Denote $B=L^{\infty}(Y)$ and given $S\subset\Lambda$, denote by $C(S)$ its centralizer in $\Lambda$.

\noindent
Suppose that $A_0\subset M^t$ is a von Neumann subalgebra such that 

$\bullet$ the inclusion $A_0\subset M^t$ has the relative property (T)

$\bullet$ $A_0\nprec_{M^t}B\rtimes\Lambda_0$, for every $\Lambda_0$ belonging to a family of subgroups $\Cal G$ of $\Lambda$.

\vskip 0.05in
\noindent Then we can find a decreasing sequence of subgroups $\{\Lambda_n\}_{n\geqslant 1}$ of $\Lambda$ with
$\Lambda_n\notin\Cal G$, for all $n\geqslant 1$,
such that $A^t\prec_{M^t}B\rtimes(\cup_{n\geqslant 1}C(\Lambda_n))$. 
\endproclaim
\vskip 0.05in 

Theorem 2 clearly follows by applying this result to the family $\Cal G$ of all amenable subgroups of $\Lambda$ in the case $t=1$ and $A_0=A$.

\vskip 0.05in
\noindent
{\it Assumptions.} (1) In order to prove Theorem 3.1 we can easily reduce to the case $t\leqslant 1$ (see e.g. the proof of Theorem 5.1). Thus, from now on, we assume that $pMp=B\rtimes\Lambda$, for some projection $p\in A$. We denote by $N:=pMp=B\rtimes\Lambda$ and by $\{v_g\}_{g\in\Lambda}\subset N$ the canonical unitaries.

\noindent (2)
We will also assume that $B\nprec_{M}A$. Indeed, otherwise by Lemma 1.3, the Cartan subalgebras $Ap$ and $B$ of $pMp$ are conjugate. Thus,  the conclusion of Theorem 3.1 automatically holds in this case.

\vskip 0.05in
Before proceeding to the proof of Theorem 3.1, let us outline it briefly in the case $p=1$. Recall from [BO08, Definition 15.1.1] that a set  $S\subset\Lambda$ is said to be {\it small relative to $\Cal G$} if $S\subset\cup_{i=1}^mg_i\Lambda_ih_i$, for some $g_i,h_i\in\Lambda$ and $\Lambda_i\in\Cal G$. We denote by $I$ the set of subsets of $\Lambda$ that are small relative to $\Cal G$.  We order $I$ by inclusion: $S\leqslant T$ iff $S\subset T$. Since $I$ is closed under finite unions, it is a directed set.
Also, we consider $\tilde M\supset M$ and the automorphisms $\{\alpha_t\}_{t\in\Bbb R}$ of $\tilde M$ constructed from the cocycle $b$ as in Section 2. 

\vskip 0.05in
\noindent {\it Outline of the proof}. The proof of Theorem 3.1 consists of two main parts:
\vskip 0.03in
\noindent
{\it Part 1}. By analyzing ``relative property (T) subsets" of $M$  we find a finite set $F\subset M$ and elements $g_S\in\Lambda\setminus S$, for every $S\in I$, such that the projection of $v_{g_S}$ onto $\sum_{x\in F}Ax$ is  uniformly  bounded away from $0$ in $||.||_2$. 

Firstly, since $A_0\nprec_{M}B\rtimes\Lambda_0$, for every $\Lambda_0\in\Cal G$, Popa's criterion  provides unitaries $a_S\in A_0$ whose support is ``almost" contained in $\Lambda\setminus S$, for every $S\in I$. Secondly, we use the fact that  $\{a_S\}_{S\in I}\subset (A_0)_1$ is a relative property (T) subset of $M$ to conclude that  for ``most" elements $g_S$ in the support of $a_S$ we have that $\alpha_t\rightarrow id$ uniformly on $\{v_{g_S}\}_{S\in I}$. Finally, since $b$ is unbounded and $B\nprec_{M}A$,  Chifan and Peterson's results  imply that $\{v_{g_S}\}_{S\in I}$ satisfy the claim.
\vskip 0.03in
\noindent
{\it Part 2}.  Let $\omega$ be a cofinal ultrafilter on $I$. We derive the conclusion by computing certain relative commutants in the ultraproduct algebra $M^{\omega}$. 

 Consider the element $g=(g_S)_S$ in the ultraproduct group $\Lambda^{\omega}$ and denote $v_g=(v_{g_S})_{S}\in M^{\omega}$. {\it Part 1} entails that the projection of $v_g$ onto $\sum_{x\in F}A^{\omega}x$ is non--zero. Let us assume for simplicity that $v_g$ in fact belongs to $A^{\omega}$.
Since $A$ is abelian, we get that $v_g$ commutes with $A$ and thus $A\subset B\rtimes\Sigma$, where $\Sigma=\Lambda\cap g\Lambda g^{-1}$. For  a set $T\subset I$, denote by $\Lambda_T$ the group generated by  $\{g_{S}g_{S'}^{-1}|S,S'\in T\}$. To reach the conclusion we combine the  following two facts: (1) an element $h\in\Lambda$ belongs to $\Sigma$ if and only if it commutes with $\Lambda_T$, for some $T\in\omega$, and (2) $\Lambda_T\notin\Cal G$, for every $T\in\omega$.
\vskip 0.05in
We are now ready to establish the first part of the proof of Theorem 3.1.

\proclaim {Lemma 3.2} In the setting of Theorem 3.1, we can find a finite set $F\subset M$ and $\delta>0$ such that the following holds: whenever $S\in I$,   there exists $g_S\in\Lambda\setminus S$ such that $\sum_{x\in F}||E_{A}(v_{g_S}x)||_2\geqslant\delta.$

\endproclaim

\noindent
{\it Remark.} 
 In the first version of this paper, we proved Theorem 3.1 and Lemma 3.2 under the  assumption that $\Gamma$ has Haagerup's property.   Stefaan Vaes pointed out to me that one can use results of [CP10] to show that Lemma 3.2 and consequently, Theorem 3.1, hold, more generally,
when $\Gamma$ has  an unbounded cocycle into a mixing representation. 

\vskip 0.05in

\noindent
{\it Proof of Lemma 3.2}.
Let $b:\Gamma\rightarrow H_{\Bbb R}$ be an unbounded cocycle. Consider $\tilde M\supset M$ and the automorphisms $\{\alpha_t\}_{t\in\Bbb R}$ of $\tilde M$ defined in Section 2.

Then the formula $\phi_t(g)=\tau(p)^{-1}\tau(\alpha_t(v_g)v_g^*)$ gives positive definite functions $\phi_t:\Lambda\rightarrow\Bbb C$. Since $||\alpha_t(v_g)-v_g||_2\rightarrow 0$, we have that   $\phi_t(g)\rightarrow 1$,  for all $g\in\Lambda$. 

 Let $\Phi_t:N\rightarrow N$ be the completely positive map defined as $\Phi_t(bv_g)=\phi_t(g)bv_g$. Then $\Phi_t$ is unital and tracial, and $||\Phi_t(x)-x||_2\rightarrow 0$, for all $x\in N$. Since the inclusion $A_0\subset N$ has the relative property (T),  for every $n\geqslant 1$ we can find $t_n>0$ such that $$||\Phi_{t_n}(a)-a||_2\leqslant \frac{||p||_2}{2^n},\hskip 0.02in\text{for all}\hskip 0.04in a\in\Cal U(A_0)\tag 3.a$$

We continue with the following:

\vskip 0.05in
\noindent
{\bf Claim.}
For any $S\in I$ and all $k\geqslant 1$, we can find $g_S\in\Lambda\setminus S$  such that $$||\alpha_{t_n}(v_{g_S})-v_{g_S}||_2\leqslant \varepsilon_n:=\sqrt{\tau(p)}\hskip 0.04in 2^{-\frac{n}{4}+2},\hskip 0.05in \forall n\in\{1,..,k\}.$$
\vskip 0.05in
\noindent
{\it Proof of the claim.} Fix $S\in I$ and $k\geqslant 1$. Then we have that $S\subset\cup_{i=1}^mg_i\Lambda_ih_i$, for some $\Lambda_i\in\Cal G$ and $g_i,h_i\in\Lambda$. Denote by $e_S$ the orthogonal projection from $L^2(N)$ onto the closed linear span of $\{Bv_g|g\in S\}$.
Since $A_0\nprec_{M}B\rtimes\Lambda_i$, for all $i$, by Remark 1.2 we can find $a_S\in\Cal U(A_0)$ with $$||e_S(a_S)||_2\leqslant\sum_{i=1}^m||E_{B\rtimes\Lambda_i}(v_{g_i}^*a_Sv_{h_i}^*)||_2\leqslant \frac{||p||_2}{2^{k}}\tag 3.b$$

Let $\tilde a_S=a_S-e_S(a_S)$. Since $||a_S||_2=||p||_2$, we get that $ ||\tilde a_S||_2>\frac{||p||_2}{2}$.
On other hand, by combining (3.a), (3.b)  and the triangle inequality we derive that $||\Phi_{t_n}(\tilde a_S)-\tilde a_S||_2\leqslant ||\Phi_{t_n}(a_S)-a_S||_2+2||e_S(a_S)||_2\leqslant 3\cdot 2^{-n}||p||_2$, for all $n\leqslant k$. We altogether deduce that  $||\Phi_{t_n}(\tilde a_S)-\tilde a_S||_2< 3\cdot 2^{-n+1}||\tilde a_S||_2$.

Now, since $\sum_{n=1}^k 2^{n-6}\cdot(3\cdot 2^{-n+1})^2=9\cdot\sum_{n=1}^k2^{-n-3}<\frac{9}{16}<1$, we get that $$\sum_{n=1}^k2^{n-6}||\Phi_{t_n}(\tilde a_S)-\tilde a_S||_2^2<||\tilde a_S||_2^2.$$

Write $\tilde a_S=\sum_{g\in\Lambda\setminus S}b_gv_g$, where $b_g\in B$. Then the last inequality rewrites as $$\sum_{g\in\Lambda\setminus S}(\sum_{n=1}^k2^{n-6}|\phi_{t_n}(g)-1|^2)\cdot||b_g||_2^2<\sum_{g\in\Lambda\setminus S}||b_g||_2^2.$$ 

Thus, 
 we can find $g_S\in\Lambda\setminus S$ satisfying $\sum_{n=1}^k2^{n-6}|\phi_{t_n}(g)-1|^2<1$. Therefore, $|\Phi_{t_n}(g)-1| < 2^{-\frac{n-6}{2}}$, for all $n\in\{1,..,k\}$. 
Finally, since $||\alpha_t(v_{g})-v_{g}||_2^2=2\tau(p)(1-\phi_t(g))$, for all $g\in\Lambda$ and $t\in\Bbb R$, the claim is proven. \hfill$\square$

\vskip 0.05in 
 Now, assume by contradiction that the conclusion of the lemma is false. Then we can find a sequence $\{S_k\}_{k\geqslant 1}\subset I$ with the following property: if  $g_k\in\Lambda\setminus S_k$, for all $k\geqslant 1$, then  
$||E_{A}(v_{g_k}x)||_2\rightarrow 0$, as $k\rightarrow\infty$, for every $x\in M$.

Let $k\geqslant 1$. By applying the above Claim to $S=S_k$ and $k$, we can find $g_{k}\in\Lambda\setminus S_k$ such that  $||\alpha_{t_n}(v_{g_{k}})-v_{g_{k}}||_2\leqslant\varepsilon_n$, for all $n\in\{1,..,k\}$. Since the map $t\rightarrow ||\alpha_t(x)-x||_2$ is a decreasing function of $|t|$, it follows that $\alpha_t\rightarrow id$ uniformly on the tail of $(v_{g_k})_{k\in\Bbb N}$.

On the other hand, as $g_k\in\Lambda\setminus S_k$, we have that $||E_{A}(v_{g_k}x)||_2\rightarrow 0$, for all $x\in M$. 
Since $v_{g_k}$ normalizes $B$, $B$ is abelian and $\alpha_t\rightarrow id$ uniformly on the tail of $(v_{g_k})_{k\in\Bbb N}$, we are in position to apply Theorem 2.4 and conclude that $\alpha_t\rightarrow id$ uniformly on $(B)_1$.	 
Since $B\nprec_{M}A$ by assumption, Theorem 2.3 gives that $\alpha_t\rightarrow id$ uniformly on $(pMp)_1$. Lemma 2.1 implies that $b$ is bounded, which provides the desired contradiction.
 \hfill$\blacksquare$
\vskip 0.05in
\noindent
{\it Remark.} Assume that $\Gamma$ has Haagerup's property, i.e. we can take the cocycle $b:\Gamma\rightarrow H_{\Bbb R}$ to be {\it proper}.
 Then Lemma 3.2 holds without assuming that  $B\nprec_{M}A$ or that $B$ is abelian. Indeed, the Claim provides $n\geqslant 1$ and $g_S\in\Lambda\setminus S$, for every $S\in I$, such that  $\inf_{S\in I}||E_M\circ\alpha_{t_n}(v_{g_S})||_2>0$. Since $b$ is proper, 
$E_M\circ\alpha_{t_n}:M\rightarrow M$ is ``compact relative to $A$". Combining these two facts readily gives the conclusion of Lemma 3.2.

As a consequence, when $\Gamma$ has Haagerup's property, Theorem 3.1 stays true if we assume that $M^t=B\rtimes\Lambda$, for an arbitrary tracial von Neumann algebra $B$.

\vskip 0.1in

\noindent{\bf 3.3 Ultraproduct algebras}.
For the second part of the proof of Theorem 3.1 we need to introduce some ultraproduct machinery (see e.g. [BO08, Appendix A]). Recall that $I$ denotes the directed set of subsets $S\subset\Lambda$ that are small relative to $\Cal G$. 

An {\it ultrafilter} $\omega$ on $I$ is  a collection of subsets of $I$ which is closed under finite unions,  does not contain the empty set and contains either $T$ or $I\setminus T$, for every subset $T$ of $I$.
Given $(x_S)_S\in\ell^{\infty}(I)$, its {\it limit along} $\omega$, denoted $\lim_{S\rightarrow\omega}x_S$, is the unique $x\in\Bbb C$ such that the set $\{S\in I|\hskip 0.02in |x_S-x|\leqslant\varepsilon\}$ belongs to $\omega$, for every $\varepsilon>0$. An ultrafilter
$\omega$ is called {\it cofinal} if it contains all the sets of the form $\{S\in I|S\geqslant S_0\}$, for some $S_0\in I$.

From now on, we fix a cofinal ultrafilter $\omega$ on $I$.
Note that $\ell^{\infty}(I,M)$ endowed with the norm $||(x_S)_S||=\sup_{S\in I}||x_S||$ is a C$^*$--algebra and that the ideal $\Cal J$ of $x=(x_S)_{S}\in\ell^{\infty}(I,M)$ satisfying $\lim_{S\rightarrow\omega}||x_S||_2=0$ is norm--closed. 
We define the {\it ultraproduct algebra}  $M^{\omega}$ as the quotient $\ell^{\infty}(I,M)/\Cal J$.
Then $M^{\omega}$ is a C$^*$--algebra and  $\tau_{\omega}:M^{\omega}\rightarrow\Bbb C$ given by $\tau_{\omega}((x_S)_S)=\lim_{S\rightarrow\omega}\tau(x_S)$ is a faithful tracial state. 

Moreover, $M^{\omega}$ is a von Neumann algebra. Indeed, the proof of [Ta03, XIV, Theorem 4.6], which deals with the particular case $I=\Bbb N$, applies verbatim for a general set $I$. Note that the trace $\tau_{\omega}$ induces a $||.||_2$ on $M^{\omega}$ given by $||(x_S)_S||_2=\lim_{S\rightarrow\omega}||x_S||_2.$
 We view $M$ as a von Neumann subalgebra of $M^{\omega}$ via the embedding $x\rightarrow (x_S)_S,$ where $x_S=x$, for all $S\in I$. 
Also, for a von Neumann subalgebra $Q$ of $M$, we view $Q^{\omega}$ as a subalgebra of $M^{\omega}$, in the natural way.

Now, recall that $N=B\rtimes\Lambda$.
 We denote by $\Lambda^{\omega}$ the ultraproduct group $(\prod_{S\in I}\Lambda)/\Cal K$, where $\Cal K=\{(g_S)_S|\lim_{S\rightarrow\omega}g_S=e\}$. If $g=(g_S)_S\in\Lambda^{\omega}$, we let
  $v_g:=(v_{g_S})_S\in\Cal U(N^{\omega})$.  Notice that this notation is consistent with the inclusion $\Lambda<\Lambda^{\omega}$.

Finally,  note that $\Lambda^{\omega}=\{v_g\}_{g\in\Lambda^{\omega}}\subset \Cal U(N^{\omega})$  normalizes $B^{\omega}$. Moreover, if $g=(g_S)_S\in\Lambda^{\omega}$, then $E_{B^{\omega}}(v_g)=(E_B(v_{g_S}))_S=(\tau(v_{g_S}))_S=\tau_{\omega}(v_g)$. Therefore, $B^{\omega}$ and $\Lambda^{\omega}$ are in a crossed product position inside $N^{\omega}$.

\vskip 0.1in
\noindent
{\it Remark}. The proof that we give below is a simplified version of our initial proof that was provided to us by Stefaan Vaes.
 \vskip 0.1in

\noindent
{\it Proof of Theorem 3.1}.  Let $g=(g_S)_S\in\Lambda^{\omega}$, where $\{g_S\}_{S\in I}$ are given by Lemma 3.2. We define $\Sigma=\Lambda\cap g\Lambda g^{-1}$ and claim that $A\prec_{M}B\rtimes\Sigma$. 

Assuming by contradiction that this is false, we can find a sequence $a_n\in\Cal U(A)$ such that $||E_{B\rtimes\Sigma}(y^*a_nx)||_2\rightarrow 0$, for any $x,y\in M$. Denote by $\Cal K\subset L^2(M^{\omega})$ the closed linear span of $Mv_gM$ and by $P$ the orthogonal projection from $L^2(M^{\omega})$ onto $\Cal K$.

Let us show that $\langle a_n\xi a_n^*,\eta\rangle\rightarrow 0$, as $n\rightarrow\infty$, for all $\xi,\eta\in\Cal K$.
To see this, it suffices to prove that $\langle a_nxv_gx'a_n^*,yv_gy'\rangle\rightarrow 0$, for all $x,x',y,y'\in M$.  Note that for every $z\in M$ we have that $E_M(v_g^*zv_g)=E_M(v_g^*E_{B\rtimes\Sigma}(z)v_g)$. Hence, we deduce that $$\langle a_nxv_gx'a_n^*,yv_gy'\rangle=\tau(v_g^*y^*a_nxv_gx'a_n^*{y'}^*)=\tau(E_M(v_g^*y^*a_nxv_g)x'a_n^*{y'}^*)=$$ $$\tau(E_M(v_g^*E_{B\rtimes\Sigma}(y^*a_nx)v_g)x'a_n^*{y'}^*).$$

Since $||E_{B\rtimes\Sigma}(y^*a_nx)||_2\rightarrow 0$, we conclude that $\langle a_nxv_gx'a_n^*,yv_gy'\rangle\rightarrow 0$, as claimed.

 Next,  Lemma 3.2 provides a finite set $F\subset M$ such that $\sum_{x\in F}||E_{A^{\omega}}(v_gx)||_2\geqslant\delta$. In particular, there is $x\in F$ such that $E_{A^{\omega}}(v_gx)\not=0$. 
 We define $\xi=P(E_{A^{\omega}}(v_gx))$ and claim that $\xi\not=0$. Since $E_{A^{\omega}}(v_gx)\not=0$, we get that $||v_gx-E_{A^{\omega}}(v_gx)||_2<||v_gx||_2$. Since $v_gx\in\Cal K$, it follows that $||v_gx-\xi||_2=||P(v_gx-E_{A^{\omega}}(v_gx))||_2<||v_gx||_2$. Hence $\xi\not=0$.

 Since $\Cal K$ is an $M$-$M$ bimodule and $A$ is abelian, we have that $a\xi=\xi a$, for all $a\in A$. In particular, we have $\langle a_n\xi a_n^*,\xi\rangle=||\xi||_2^2,$ for all $n$. This contradicts the fact that $\langle a_n\xi a_n^*,\xi\rangle\rightarrow 0$ and proves that $A\prec_{M}B\rtimes\Sigma$.

To finish the proof it suffices to produce a decreasing sequence $\{\Lambda_n\}_{n\geqslant 1}$ of subgroups of $\Lambda$ such that $\Lambda_n\notin\Cal G$, for all $n\geqslant 1$, and $\Sigma=\cup_{n\geqslant 1}C(\Lambda_n)$.

Next, for $T\subset I$, we let $\Lambda_T$ be the subgroup of $\Lambda$ generated by $\{g_{S}g_{S'}^{-1}|S,S'\in T\}$.
It is clear that an element $h\in\Lambda$ belongs to $\Sigma$ if and only if there exists $T\in\omega$ such that $h\in C(\Lambda_T)$. Thus, if we enumerate $\Sigma=\{h_n\}_{n\geqslant 1}$, then for every $n\geqslant 1$ there exists $T_n\in\omega$ such that $h_n\in C(\Lambda_{T_n})$. Put $W_n=\cap_{i=1}^nT_i$. Then $W_n\in\omega$ and $W_n\supset W_{n+1}$ for all $n\geqslant 1$, and we have that $\Sigma=\cup_{n\geqslant 1}C(\Lambda_{W_n})$.

Finally, let us argue that $\Lambda_W\notin\Cal G$, for every $W\in\omega$. Assume by contradiction that $\Lambda_W\in\Cal G$ and fix $S'\in W$. Then 
the set $S''=\Lambda_Wg_{S'}$ is small relative to $\Cal G$, i.e. $S''\in I$. Since $\omega$ is a cofinal ultrafilter on $I$ and $W\in\omega$, we can find $S\in W$ such that $S\supset S''$. Since $g_S\in\Lambda_Wg_{S'}=S''$ this contradicts the fact that $g_S\in\Lambda\setminus S$.
\hfill$\blacksquare$
\vskip 0.1in

Next, we notice that the proof of Theorem 3.1 also yields the following:

\proclaim {Lemma 3.4} Let $(B,\tau)$ be a tracial von Neumann algebra and $\Lambda\curvearrowright B$ be a trace preserving action. Let $N=B\rtimes\Lambda$ and $A\subset N$ be an abelian von Neumann subalgebra. Assume that we can find two sequences $\{a_n\}_{n\geqslant 1}\subset (A)_1$ and $\{g_n\}_{n\geqslant 1}\subset\Lambda$ such that  $g_n\rightarrow\infty$ and $\inf_{n}||E_B(a_nv_{g_n}^*)||_2>0$. 

\vskip 0.03in
\noindent
Then we can find  a decreasing sequence $\{\Lambda_n\}_{n\geqslant 1}$ of infinite subgroups of $\Lambda$ such that $A\prec_{N}B\rtimes(\cup_{n\geqslant 1}C(\Lambda_n))$.
\endproclaim

\noindent
{\it Proof.} Let $\omega$ be a free ultrafilter on $\Bbb N$ and consider the notations from 3.3 for $I=\Bbb N$ . Put  $g=(g_n)_n\in\Lambda^{\omega}$. The hypothesis guarantees that $b:=E_{B^{\omega}}(av_g^*)\not=0$.  This implies that $E_{A^{\omega}}(bv_g)\not=0$. 

Let $\Sigma=\Lambda\cap g\Lambda g^{-1}$. We claim that $A\prec_{M}B\rtimes\Sigma$. 
The claim follows by adjusting the proof of Theorem 3.1. Assuming by contradiction that the claim is false we can find $a_n\in\Cal U(M)$ such that $||E_{B\rtimes\Sigma}(y^*a_nx)||_2\rightarrow 0$, for all $x,y\in M$.
  Let $x,x',y,y'\in M$. Since $E_M(v_g^*b^*zbv_g)=E_M(v_g^*b^*E_{B\rtimes\Sigma}(z)bv_g)$, for every $z\in M$, we deduce that $$|\langle a_n(xbv_gx')a_n^*,ybv_gy'\rangle|=|\tau(v_g^*b^*y^*a_nxbv_gx'a_n^*{y'}^*)|=$$ $$|\tau(E_M(v_g^*b^*y^*a_nxbv_g)x'a_n^*{y'}^*)|=|\tau(E_M(v_g^*b^*E_{B\rtimes\Sigma}(y^*a_nx)bv_g)x'a_n^*{y'}^*)|\rightarrow 0.$$
  
  Denote by $\Cal K\subset L^2(M^{\omega})$ the closed linear span of $Mbv_gM$. The above calculation shows that
 $\langle a_n\xi a_n^*,\eta\rangle\rightarrow 0$, for all $\xi,\eta\in\Cal K$.
  By the proof of Theorem 3.1, this is enough to imply that $A\prec_{M}B\rtimes\Sigma$.
 
The proof of Theorem 3.1 also gives that $\Sigma=\cup_{n\geqslant 1}C(\Lambda_{W_n})$, for some decreasing sequence $\{W_n\}_{n\geqslant 1}$ of sets $W_n\in\omega$. Since every set 
in $\omega$ is infinite, it follows that $\Lambda_{W_n}$ is infinite, for all $n$.
\hfill$\blacksquare$

\vskip 0.05in

We end this section with a consequence of Theorem 3.1 and a result of Ozawa [Oz08]. We say that a group $\Lambda$ has {\it Haagerup's property relative to a subgroup} $\Sigma$ if we can find a sequence $\phi_n:\Lambda\rightarrow\Bbb C$ of positive definite functions such that 

$\bullet$ for all $g\in\Lambda$, we have that $\phi_n(g)\rightarrow 1$,  and

$\bullet$ for all $n\geqslant 1$ and $\varepsilon>0$, we can find $g_1,..,g_k,h_1,..,h_k\in \Lambda$ such that $|\phi_n(g)|<\varepsilon$, 

\hskip 0.12in for all $g\in\Lambda\setminus (\cup_{i=1}^kg_i\Sigma h_i)$.

\proclaim {Corollary 3.5}  Let $\Gamma<SL_2(\Bbb Z)$ be a non--amenable subgroup. Denote $M=L(\Bbb Z^2\rtimes\Gamma)$. Let $\Lambda$ be a countable group such that $M=L\Lambda$.

\vskip 0.02in
\noindent 
Then $\Lambda$ has Haagerup's property relative to some infinite amenable subgroup $\Sigma$.

\endproclaim

\noindent
{\it Proof.}  Since the inclusion $L(\Bbb Z^2)\subset M$ has the relative property (T) ([Bu91],[Po01]) and $\Gamma$ has Haagerup's property, by the remark just before subsection 3.3 we are in position to apply Theorem 3.1. By applying Theorem 3.1 in the case $B=\Bbb C1$ and $\Cal G$ is the family of finite subgroups of $\Lambda$  we get that $L(\Bbb Z^2)\prec_{M} L(\Sigma)$, where $\Sigma=\cup_{n\geqslant 1}C(\Lambda_n)$, for some decreasing sequence $\{\Lambda_n\}_{n\geqslant 1}$ of infinite subgroups of $\Lambda$.
On the other hand, by [Oz08] we have that $M$ is solid, i.e. the commutant of any diffuse subalgebra is amenable. It follows that $C(\Lambda_n)$ is amenable, for all $n\geqslant 1$, and thus $\Sigma$ is amenable.

Now, since $L(\Bbb Z^2)\subset M$ is a Cartan subalgebra and $L(\Bbb Z^2)\prec_{M}L(\Sigma)$, we can find $x_1,..,x_n,y_1,..,y_n\in M$ such that $(L(\Bbb Z^2))_1$ is contained in the linear span of $\{x_i(L(\Sigma))_1y_i|\hskip 0.03ini\in\{1,..,n\}\}.$ By using again that $\Gamma$ has Haagerup's property, the conclusion follows easily.
\hfill$\blacksquare$\vskip 0.2in

\head \S 4.  A conjugacy criterion for Cartan subalgebras.\endhead
\vskip 0.2in

In this section we prove a general criterion for unitary conjugacy of Cartan subalgebras and derive Theorem 3 as a corollary. 

Before stating our criterion, let us recall from  [Ga02, Definition I.5]  the notion of  cost of an equivalence relation.
Let $\Cal R$ be a countable, measure preserving equivalence relation on a standard probability space $(X,\mu)$. A countable family $\Theta=\{\theta_i:Y_i\rightarrow Z_i\}_{i\in I}\subset [[\Cal R]]$ is a {\it graphing of}  $\Cal R$, if $\Cal R$ is the smallest equivalence relation $\Cal S$ satisfying $\theta_i\in[[\Cal S]]$, for all $i\in I$. The cost of a graphing $\Theta$ is defined as $\Cal C(\Theta)=\sum_{i\in I}\mu(Y_i)$. Finally, the {\it cost} of $\Cal R$ is defined by 
$\Cal C(\Cal R)=\inf\{\Cal C(\Theta)|\hskip 0.02in \Theta$ is a graphing of $\Cal R\}$.

\proclaim {Theorem 4.1} Let $A$  be a Cartan subalgebra of a separable II$_1$ factor $M$. 
Assume that the equivalence relation $\Cal R$  associated with the inclusion $(A\subset M)$ satisfies $\Cal C(\Cal R)>1$.

\noindent
Let $B\subset M$ be a Cartan subalgebra. 
 Suppose that there is an amenable von Neumann subalgebra $N\subset M$ such that either

\noindent
(1) $A\subset N$ and $B\prec_{M}N$, or

\noindent
(2) $A\prec_{M}N$ and $B\subset N$.

\vskip 0.02in
\noindent
Then we can find a unitary element $u\in M$ such that $uAu^*=B$.
\endproclaim

Before proceeding to the proof of Theorem 4.1 let us derive Theorem 3 from it. We moreover prove a generalization of Theorem 3 which involves amplifications.

\proclaim {Theorem 4.2} Let $\Gamma\curvearrowright (X,\mu)$ be a free ergodic p.m.p. action and assume that $\beta_1^{(2)}(\Gamma)>0$. Denote   $A=L^{\infty}(X)$ and $M=L^{\infty}(X)\rtimes\Gamma$.   Let  $B\subset M^t$ be a Cartan subalgebra, for some $t>0$.

\vskip 0.02in
\noindent
 If  there exists an amenable von Neumann subalgebra $N$ of $M^t$ such that $A^t\prec_{M^t}N$ and $B\subset N$,
then we can find a unitary element $u\in M^t$ such that $uA^tu^*=B$.
\endproclaim

\noindent
{\it Proof.}  Let $\Cal R$ be the equivalence relation induced by the action $\Gamma\curvearrowright X$.  Then [Ga01, Corollaire 3.23 and Corollaire 3.16] give that $\Cal C(\Cal R)\geqslant\beta_1^{(2)}(\Cal R)+1=\beta_{1}^{(2)}(\Gamma)+1$ and thus $\Cal C(\Cal R)>1$. 
This inequality and [Ga99, Proposition II.6] imply that $\Cal C(\Cal R^t)>1$, for every $t>0$. Since $\Cal R^t$ is precisely the equivalence relation of the inclusion $(L^{\infty}(X)^t\subset M^t)$,  the conclusion follows by applying Theorem 4.1. \hfill$\blacksquare$

\vskip 0.05in
\noindent
As a first step towards Theorem 4.1 we show that conditions (1) and (2) are equivalent.

\proclaim {Proposition 4.3}  If $A$ and $B$ are Cartan subalgebras of a separable II$_1$ factor $M$, then the following are equivalent:

\noindent
(1) there is an amenable subalgebra $N\subset M$ such that $A\subset N$ and 
 $B\prec_{M}N$.

\noindent
(2) there is an amenable subalgebra $N\subset M$  such that $A\prec_{M} N$ and $B\subset N$.

\noindent
(3) there is  an amenable subalgebra $N\subset rMr$, for some non--zero projection $r\in M$, such that $A\prec_{M}Ns$ and $B\prec_{M}Ns$, for every non--zero projection $s\in N'\cap rMr$.
\endproclaim

\noindent 
{\it Proof.} By symmetry, it suffices to show that  (1) implies (3) and that (3) implies (1).
\vskip 0.05in

\noindent 
(1) $\Longrightarrow$ (3). Let $N\subset M$ amenable such that  $A\subset N$ and $B\prec_{M} N$. By a maximality argument, we can find a non--zero projection $r \in N'\cap M$ such that $B\prec_{M}Ns$, for any  non--zero projection $s\in N'\cap M$ with $s\leqslant r$. Since $A\subset N$, we also have that $A\prec_{M}Ns$, for every non--zero projection $s\in N'\cap M$. It follows that (3) holds for $Nr\subset rMr$.

\vskip 0.05in
\noindent (3) $\Longrightarrow$ (1). Let $N\subset rMr$ satisfying (3).
Since $A\prec_{M}N$, we can find  projections $p\in A,q\in N$, a $*$--homomorphism $\psi:Ap\rightarrow qNq$ and a non--zero partial isometry $v\in qMp$ such that $\psi(x)v=vx$, for all $x\in Ap$, $v^*v=p$ and $q':=vv^*\in\psi(Ap)'\cap qMq$. Moreover, by Lemma 1.5 we may assume that $\psi(Ap)$ is maximal abelian in $qNq$.

Let $P$ be the von Neumann algebra generated by the normalizer of $\psi(Ap)$ in $qNq$. 
Also, let $Q\subset pMp$ be the von Neumann algebra generated by $v^*Pv$. We have that
\vskip 0.05in
\noindent
{\bf Claim 1.} $B\prec_{M}Q$.
\vskip 0.05in
\noindent
{\bf Claim 2.} $Q$ is amenable.
\vskip 0.05in
Before proving these claims let us indicate how they imply the conclusion. Firstly, since $v^*\psi(Ap)v=Ap$,  we have that $Ap\subset Q$. 
Since $Q$ is amenable and $Ap\subset Q$, we can construct an amenable subalgebra $R\subset M$ such that $A\subset R$, $p\in R$ and $pRp=Q$. Since $B\prec_{M}Q$, it follows that $B\prec_{M}R$ and therefore (1) holds.

\vskip 0.05in
\noindent
{\it Proof of Claim 1.} By 
Lemma 1.6 (2) we deduce that $P\prec_{M}Q$. By a maximality argument we
 can find a non--zero projection $e\in P'\cap qNq$ such that $Pf\prec_{M}Q$, for any non--zero projection $f\in P'\cap qNq$ satisfying $f\leqslant e$.

Next, for $u\in\Cal N_{pMp}(Ap)$, define $\theta_u\in$ Aut$(Ap)$ by $\theta_u(x)=uxu^*$. Then for any $y\in\psi(Ap)$ we have that $vuv^*y=(\psi\circ\theta_u\circ\psi^{-1})(y)vuv^*$.
Since $\psi(Ap)$ is maximal abelian in $qNq$, it follows that $E_N(vuv^*)\in P$. Since $Ap$ is regular in $pMp$, we get that $E_N(q'Mq')\subset P$. 
Since $e\in P'\cap qNq$, Lemma 1.6 (1) gives that $N\prec_{N}Pe$.  
 By  [Va07, Lemma 3.7], the combination of the last two paragraphs implies that $N\prec_{M}Q$.

Thus, we can find a non--zero projection $s\in N'\cap rMr$ such that $Nt\prec_{M}Q$, for every non--zero projection $t\in N'\cap rMr$ with $t\leqslant s$.
Since $B\prec_{M} Ns$, by our assumption,  applying [Va07, Lemma 3.7] again yields that $B\prec_{M}Q$.\hfill$\square$

\vskip 0.05in
\noindent
{\it Proof of Claim 2}.
We start by identifying $Ap=L^{\infty}(T)$ and $\psi(Ap)=L^{\infty}(W)$, where $T,W$ are probability spaces. Let $\theta:W\rightarrow T$ be a probability space isomorphism such that $\psi(x)=x\circ\theta$, for all $x\in Ap=L^{\infty}(T).$  Let $\Cal R$ be the equivalence relation on $W$ associated with the Cartan subalgebra inclusion $(\psi(Ap)\subset P)$ ([FM77]). 
Since $N$ and hence $P$ is amenable, we get that $\Cal R$ is hyperfinite ([CFW81]).

Now, let $\Cal S$ be the equivalence relation on $T$ associated with the inclusion $Ap\subset pMp$. Set $\Cal S_0=\Cal S\cap (\theta\times\theta)(\Cal R)$. Then $\Cal S_0$ is a hyperfinite subequivalence relation of $\Cal S$. By [FM77, Theorem 1], we can find an amenable von Neumann subalgebra $Q_0\subset pMp$ such that $Ap\subset Q_0$ and $\Cal S_0$ is the equivalence relation associated to the inclusion $Ap\subset Q_0$.

We claim that $Q\subset Q_0$, which implies that $Q$ is amenable. Let $u\in\Cal N_{qNq}(\psi(Ap))$ and define $\phi\in [\Cal R]$ by $y\circ\phi=uyu^*$, for all $y\in\psi(Ap)$. Denote $\alpha=\theta\phi\theta^{-1}\in$ Aut$(T)$ and $w=v^*uv$.
Then  we have $wx=(x\circ\alpha)w$, for  every $x\in Ap$.

Since $Ap\subset pMp$ is maximal abelian, the left and right supports of $w$ lie in $Ap$. Thus, $ww^*=1_{T_1},w^*w=1_{T_2}$, where $T_1,T_2\subset T$ are Borel. Then $\alpha(T_1)=T_2$ and $\beta:=\alpha_{|T_1}$ belongs to $[[\Cal S]]$. Moreover, $w\in Au_{\beta}^*$, where $u_{\beta}\in pMp$ is the partial isometry implementing $\beta$. Finally, since $\beta$ belongs to $\theta [[\Cal R]]\theta^{-1}\cap [[\Cal S]]=[[\Cal S_0]]$, we get that $u_{\beta}\in Q_0$. Thus, $w=v^*uv\in Q_0$, for all $u\in\Cal N_{qNq}(\psi(Ap))$ and hence $Q\subset Q_0$.
\hfill$\blacksquare$

\vskip 0.05in

Next, we introduce a notion of quasi--normality for subequivalence relations which is inspired by Popa's notion of {\it wq--normal} subgroups ([Po04, Definition 2.3]) and by Peterson and Thom's notion of {\it s--normal} subgroupoids ([PT07, Definition 6.3]).
\vskip 0.05in
\noindent{\it Definition 4.4} Let $\Cal S\subset \Cal R$ be countable measure preserving equivalence relations on a probability space $(X,\mu)$.  We say that $\Cal S$ is  {\it q--normal} in $\Cal R$ if we can find $\theta_n\in [[\Cal R]]$, with $\theta_n:Y_n\rightarrow Z_n$, for all $n\geqslant 1$, such that

\vskip 0.02in
(1) $\{\theta_n\}_{n\geqslant 1}$ generate $\Cal R$ as an equivalence relation and

(2) the equivalence relation \hskip 0.05in $\{(x,y)\in Y_n\times Y_n|\hskip 0.02in (x,y)\in\Cal S$\hskip 0.05in and\hskip 0.05in $(\theta_n(x),\theta_n(y))\in\Cal S\}$ 

\hskip 0.23in has infinite orbits, for all $n\geqslant 1$.
\vskip 0.05in

We continue with a result which will be essential in the proof of Theorem 4.1.

\proclaim {Proposition 4.5} Let $M$ be a separable II$_1$ factor together with two  Cartan subalgebras $A$ and $B$. Suppose that there is no unitary $u\in M$ such that $uAu^*=B$.
Assume that there is an amenable von Neumann subalgebra $N\subset M$ such that  $A\subset N$ and $B\prec_{M}N$. 

\noindent
Identify $A=L^{\infty}(X)$, where $(X,\mu)$ is a probability space. Denote by $\Cal R$ and $\Cal S$ the equivalence relations on $X$ associated with the inclusions $A\subset M$ and  $A\subset N$.
\vskip 0.03in

\noindent
Then we can find a set $X_0\subset X$ of positive measure, an equivalence relation $\Cal T$ on $X_0$ with $\Cal S_{|X_0}\subset\Cal T\subset \Cal R_{|X_0}$ and a partition $\{X_k\}_{k\geqslant 1}$ of $X_0$ into Borel subsets such that  

\vskip 0.02in
\noindent
(1) $\Cal S_{|X_0}$ is hyperfinite and its restriction to any Borel set of positive measure has infinite orbits,

\noindent
(2) $\Cal S_{|X_0}$ is q--normal in $\Cal T$, and 

\noindent
(3) almost every $\Cal R_{|X_k}$--class contains only finitely many $\Cal T_{|X_k}$--classes, for all $k\geqslant 1$.
\endproclaim

\noindent
{\it Proof}. Let  $N\subset M$ amenable such that $A\subset N$ and $B\prec_{M}N$. Since $A$ and $B$ are not conjugate by a unitary, by Lemma 1.3 we have that $B\nprec_{M}A$.
Then we can find projections $p\in B,q\in N$, a $*$--homomorphism $\psi:Bp\rightarrow qNq$ and a non--zero partial isometry $v\in qMp$ such that $v^*v=p$ and $\psi(b)v=vb$, for all $b\in Bp$. Since $B\nprec_{M}A$, we may also assume that $\psi(Bp)\nprec_{M}A$ ([Va07, Remark 3.8.]). Let $q'=vv^*\leqslant q$. 
\vskip 0.05in
\noindent
Before continuing we need to introduce some notations:

$\bullet$ Denote by $P$ the von Neumann algebra generated by $A$ and $q'Mq'$.

$\bullet$ Denote by $\Cal R_0$ the equivalence relations on $X$ associated with the inclusion $A\subset P$.

$\bullet$ For  $\phi\in [[\Cal R]]$, let $u_{\phi}\in M$ be a partial isometry which implements $\phi$.

$\bullet$ Fix a sequence $\{\phi_m\}_{m\geqslant 1}\subset [[\Cal R_0]]$ such that $\Cal R_0=\sqcup_{m\geqslant 1}\{(\phi_m(x),x)|x\in X\}$.

$\bullet$ Fix a sequence $\{u_n\}_{n\geqslant 1}\subset\Cal N_{pMp}(Bp)$ which generates $pMp$ as a von Neumann 

\hskip 0.13in algebra (such a sequence exists because $Bp$ is regular in $pMp$). 

\vskip 0.05in
The choice of $\{\phi_m\}_{m\geqslant 1}$ guarantees that  $\{u_{\phi_m}\}_{m\geqslant 1}$ is an orthonormal basis for $P$ over $A$ (see e.g. [PP86]). Since $vu_nv^*\in q'Mq'\subset P$, we have that $vu_nv^*=\sum_{m\geqslant 1}a_{m,n}u_{\phi_m}$, where  $a_{m,n}=E_A(vu_nv^*u_{\phi_m}^*)$ and the sum converges in $||.||_2$. Let $X_{m,n}\subset X$ be the essential support of $a_{m,n}$ and $\phi_{m,n}$ be the restriction of $\phi_m$ to $\phi_m^{-1}(X_{m,n})$. Hence, there is a partial isometry $v_{m,n}\in A$ with support $X_{m,n}$ such that $1_{X_{m,n}}u_{\phi_m}=v_{m,n}u_{\phi_{m,n}}$.
Altogether, we get that $vu_nv^*=\sum_{m\geqslant 1}a_{m,n}v_{m,n}u_{\phi_{m,n}},$ for all $n\geqslant 1$.

Since $q'Mq'=v(pMp)v^*$, we have that $P$ is generated by $A$ and $\{vu_nv^*\}_{n\geqslant 1}$. The last identity in the previous paragraph implies that $P$ is generated by $A$ and $u_{\phi_{m,n}}$. We deduce that $\Cal R_0$ is generated, as an equivalence relation, by $\{\phi_{m,n}\}_{m,n\geqslant 1}$ and id$_X$.
\vskip 0.05in
The proof is divided between  three claims. 
The first and most important claim asserts that each $\phi_{m,n}$ ``quasi--normalizes" $\Cal S$.
\vskip 0.05in
\noindent
{\bf Claim 1.} Fix $m,n\geqslant 1$. Let $Y$ be the domain of $\phi_{m,n}$. Then the equivalence relation 
$\{(x,y)\in Y\times Y|\hskip 0.02in(x,y)\in\Cal S\hskip 0.05in\text{and}\hskip 0.05in(\phi_{m,n}(x),\phi_{m,n}(y))\in\Cal S\}$
 has infinite orbits.

\vskip 0.05in
\noindent
{\it Proof of claim 1}. Assume by contradiction that the claim is false. Then we can find a Borel set $Z\subset Y$ with $\mu(Z)>0$ such that $\phi={\phi_{m,n}}_{|Z}$ satisfies 
$(\phi(x),\phi(y))\notin\Cal S$, for all $(x,y)\in\Cal S\cap (Z\times Z)$ with $x\not=y$. 

Let us show that there is $a\in A$ such that $\delta=\langle au_{\phi},vu_nv^*\rangle>0$. Since $\phi={\phi_{m}}_{|Z}$ we can find a partial isometry $c\in A$ with support $\phi_m(Z)$ such that $u_{\phi}=cu_{\phi_{m}}$.
As the projection of $vu_nv^*$ onto the closure of $Au_{\phi_m}$ is equal to $a_{m,n}u_{\phi_m}$, the projection of $vu_nv^*$ onto the closure of $Au_{{\phi_{m}}_{|Z}}$ is equal to $1_{\phi_m(Z)}a_{m,n}u_{\phi_m}=c^*a_{m,n}u_{\phi}$. Since $\phi_m(Z)$ is contained in the support of $a_{m,n}$, the latter is non--zero. Thus, $a=c^*a_{m,n}\in A$ works.  

Now, fix $b\in\Cal U(\psi(Bp))$ and set $\rho=\psi$ $\circ$ Ad$(u_n)\circ\psi^{-1}\in$ Aut$(\psi(Bp))$.  Then we have that $\rho(b)(vu_nv^*)=(vu_nv^*)b$. Since $b\in\Cal U(qMq)$ and $vu_nv^*\in qMq$,
  we have that $$\Re\hskip 0.02in\langle au_{\phi}b,\rho(b)vu_nv^*\rangle=\Re\hskip 0.02in\langle au_{\phi}b,vu_nv^*b\rangle=\Re\hskip 0.02in\langle au_{\phi},vu_nv^*\rangle=\delta>0\tag 4.a$$

On the other hand, since $a,\rho(b)\in N$ and we have that

$$\Re\hskip 0.02in\langle au_{\phi}b,\rho(b)vu_nv^*\rangle=\Re\hskip 0.02in\tau(\rho(b)^*au_{\phi}bvu_n^*v^*)\leqslant ||a||_2\hskip 0.02in||E_N(u_{\phi}bvu_n^*v^*)||_2\tag 4.b$$
By combining (4.a) and (4.b) we get that $$||E_N(u_{\phi}bvu_n^*v^*)||_2\geqslant \frac{\delta}{||a||_2},\hskip 0.05in\forall b\in\Cal U(\psi(Bp))\tag 4.c$$

\noindent
Since $\psi(Bp)\nprec_{M}A$, by Theorem 1.1 we can find a sequence $b_k\in\Cal U(\psi(Bp))$ such that $||E_A(b_kw)||_2\rightarrow 0$, for every $w\in M$. 
Let us show that $$||E_N(u_{\phi}b_kz)||_2\rightarrow 0,\hskip 0.05in\forall z\in M\tag 4.d$$ It is clear that (4.d)  contradicts (4.c) and therefore proves the claim. By Kaplansky's density theorem it is enough to prove (4.d) when $z=u_{\phi'}$, for some $\phi'\in [\Cal R]$.

Let $\{\alpha_l\}_{l\geqslant 1}\subset [[\Cal S]]$ be a sequence such that $\{u_{\alpha_l}\}_{l\geqslant 1}$ is an orthonormal basis for $N$ over $A$.  Let $X_l$ be the set of $x\in X$ for which $\phi\alpha_l\phi'(x)$ is defined and $(\phi\alpha_l\phi'(x),x)\in\Cal S$. 
We have that the sets $\{X_l\}_{l\geqslant 1}$ are mutually disjoint. Indeed, if $x\in X_{l}\cap X_{l'}$, then $(\phi(\alpha_l\phi'(x)),\phi(\alpha_{l'}\phi'(x))\in\Cal S$. Since $\alpha_l,\alpha_l'\in [[\Cal S]]$ we also have that $(\alpha_l\phi'(x),\alpha_{l'}\phi'(x))\in\Cal S$. Thus, we deduce that $\alpha_l\phi'(x)=\alpha_{l'}\phi'(x)$, hence $l=l'$.

Let $\varepsilon>0$ and $L\geqslant 1$ such that $\sum_{l\geqslant L}\mu(X_l)\leqslant\varepsilon$. Since $b_k\in \psi(Bp)\subset N$, we can write $b_k=\sum_{l\geqslant 1}E_A(b_ku_{\alpha_l}^*)u_{\alpha_l}$ and thus $E_N(u_{\phi}b_ku_{\phi'})=\sum_{l\geqslant 1}\phi(E_A(b_ku_{\alpha_l}^*))E_N(u_{\phi\alpha_l\phi'})$.
Further, since $||E_A(b_ku_{\alpha_l}^*)||\leqslant 1$ and $E_N(u_{\phi\alpha_l\phi'})=1_{X_l}u_{\phi\alpha_l\phi'}$, it follows that for all $k\geqslant 1$ we have that $$||E_N(u_{\phi}b_ku_{\phi'})||_2^2=\sum_{l\geqslant 1}||1_{X_l}\phi(E_A(b_ku_{\alpha_l}^*))||_2^2\leqslant$$ $$\sum_{l\geqslant L}||1_{X_l}||_2^2+ \sum_{l<L}||E_A(b_ku_{\alpha_l}^*)||_2^2\leqslant\varepsilon+\sum_{l<L}||E_A(b_ku_{\alpha_l}^*)||_2^2.$$
As $||E_A(b_ku_{\alpha_l}^*)||_2\rightarrow 0$, for all $l\geqslant 1$, we get that $\limsup_{k\rightarrow\infty}||E_N(u_{\phi}b_ku_{\phi'})||_2\leqslant\sqrt{\varepsilon}$.  Since $\varepsilon>0$ is arbitrary, we conclude that $||E_N(u_{\phi}b_ku_{\phi'})||_2\rightarrow 0$.\hfill$\square$

\vskip 0.05in

Next, let $q_0$ be the support projection of $E_A(q')$. Write  $q_0=1_{X_0}$, for $X_0\subset X$ Borel.

\vskip 0.05in
\noindent{\bf Claim 2.}  We can find a partition $\{X_k\}_{k\geqslant 1}$ of $X_0$ into Borel sets such that almost every $\Cal R_{|X_k}$--class contains only finitely many ${\Cal R_0}_{|X_k}$--classes, for all $k\geqslant 1$. 
\vskip 0.05in
\noindent
{\it Proof of Claim 2.}
By using a maximality argument, it suffices to prove that whenever $X_1\subset X_0$ is a set of positive measure, we can find a set $X_2\subset X_1$ of positive measure such that every $\Cal R_{|X_2}$--class contains only finitely many ${\Cal R_0}_{|X_2}$--classes.

 To see this, put $q_1=1_{X_1}$.
Since $P$ contains $q'Mq'$, we get that $q_1Pq_1$ contains $q_1q'Mq'q_1$. Thus, if $q_2$ denotes the left support of $q'q_1$, then $q_1Pq_1$ contains $w(q_2Mq_2)w^*$, for some unitary element $w\in M$. Since $q'q_1\not=0$, we have $q_2\not=0$, and it follows that $M\prec_{M}q_1Pq_1$.
Thus,  $M\prec_{M}\tilde P=q_1Pq_1\oplus A(1-q_1)$.  
 Now, the equivalence relation of the inclusion $A\subset\tilde P$ is equal to ${\Cal R_0}_{|X_1}\cup$ id$_{X\setminus X_1}$. By applying Lemma 1.8 (to the case $N=M$) our claim follows.\hfill$\square$
\vskip 0.05in

\noindent
{\bf Claim 3.} $\Cal S_{|X_0}$ is hyperfinite and its restriction to any Borel set  of positive measure has infinite orbits.  
\vskip 0.05in
\noindent
{\it Proof of Claim 3.} 
 Since $\Cal S_{|X_0}$ is the equivalence relation of the inclusion 
$(Aq_0\subset q_0Nq_0)$ and $N$ is amenable, by [CFW81] we deduce that $\Cal S_{|X_0}$ is hyperfinite.

Now, let $Y\subset X_0$ be a set of positive measure and set $r=1_{Y}$. In order to show that $\Cal S_{|Y}$ has infinite orbits it suffices to argue that $rNr\nprec_{N}A$.

 Since $\psi(Bp)\nprec_{N}A$, we get that $qNq\nprec_{N}A$. It follows that $Nq_1\nprec_{N}A$, where $q_1$ is the central support of $q$ in $N$. If $\Cal Z$ denotes the center of $N$, then $q_1$ is precisely the support of $E_{\Cal Z}(q)$. Let $q_2$ be the support of $E_A(q)$. Since $\Cal Z\subset A$, we have that $q_2\leqslant q_1$. Also, since $q'\leqslant q$ and $q_0$ is the support of $E_A(q')$, we get that $q_0\leqslant q_2$. Altogether, we derive that $q_0\leqslant q_1$. Thus, $q_0Nq_0\nprec_{N}A$ and since $r\leqslant q_0$, we get that $rNr\nprec_{N}A$. \hfill$\square$
\vskip 0.05in

We are now ready to combine all the claims and finish the proof of Proposition 4.5.
Let $\Cal T$ be the equivalence relation on $X_0$ generated by $\Cal S_{|X_0}$ and ${\Cal R_0}_{|X_0}$. Since the domain and image of each $\phi_{m,n}$ is contained in $X_0$, we get that $\Cal T$ is generated by $\Cal S_{|X_0}$ and $\{\phi_{m,n}\}_{m,n\geqslant 1}$. Since $\Cal S_{|X_0}$ has infinite orbits, Claim 1 implies that the inclusion $\Cal S_{|X_0}\subset\Cal T$ is q--normal, hence condition (2) of the conclusion is verified. Since conditions (1) and (3) also hold by claims 3 and 2, we are done.\hfill$\blacksquare$
\vskip 0.05in

The last ingredient in the proof of Theorem 4.1. is a lemma due to D. Gaboriau which asserts that cost does not increase by passing to q--normal extensions.

\proclaim {Lemma 4.6 [Ga99, Lemma V.3.]} Let $\Cal R$ be a countable, measure preserving equivalence relation on a probability space $(X,\mu)$. If $\Cal S\subset\Cal R$ is a q--normal subequivalence relation, then $\Cal C(\Cal R)\leqslant\Cal C(\Cal S)$. 
\endproclaim
\noindent
{\it Proof.} For the reader's convenience let us recall from [Ga99] the proof of this lemma. Let $\varepsilon>0$ and $\Theta$ be a graphing of $\Cal S$ such that $\Cal C(\Theta)\leqslant \Cal C(\Cal S)+\frac{\varepsilon}{2}$. Since $\Cal S$ is q--normal in $\Cal R$, we can find a sequence $\{\theta_n:Y_n\rightarrow Z_n\}_{n\geqslant 1}\subset [[\Cal R]]$ which generates $\Cal R$ as an equivalence relation such that $\Cal S_n=\{(x,y)\in (Y_n\times Y_n)\cap\Cal S|\hskip 0.02in (\theta_n(x),\theta_n(y))\in\Cal S\}$ has infinite orbits, for all $n\geqslant 1$. 
Let $Y_n^0\subset Y_n$ be a Borel set of measure at most $\frac{\varepsilon}{2^{n+1}}$ that intersects almost every $\Cal S_n$--class. 

We claim that $\tilde{\Theta}=\Theta\cup\{{\theta_n}_{|Y_n^0}\}_{n\geqslant 1}$ is a graphing for $\Cal R$. Let $\Cal R_0\subset\Cal R$
be the equivalence relation generated by $\tilde{\Theta}$. For $n\geqslant 1$ and almost every $x\in Y_n$ we can find $y\in Y_n^0$ such that $(x,y)\in\Cal S_n$. Since $\Cal S\subset\Cal R_0$, we get that $(x,y),(\theta_n(x),\theta_n(y))\in\Cal R_0$. Also, since ${\theta_n}_{|Y_n^0}\in[[\Cal R_0]]$, we have that $(y,\theta_n(y))\in\Cal R_0$. Altogether, it follows that $(x,\theta_n(x))\in\Cal R_0$. Since $\{\theta_n\}_{n\geqslant 1}$ generates $\Cal R$, we deduce that $\Cal R_0=\Cal R$, as claimed. 	

Now, $\Cal C(\tilde\Theta)=\Cal C(\Theta)+\sum_{n\geqslant 1}\mu(Y_n^0)\leqslant\Cal C(\Theta)+\frac{\varepsilon}{2}\leqslant \Cal C(\Cal S)+\varepsilon$. Since $\tilde{\Theta}$ is a graphing for $\Cal R$, we get that $\Cal  C(\Cal R)\leqslant \Cal C(\tilde{\Theta})\leqslant \Cal C(\Cal S)+\varepsilon$. As $\varepsilon>0$ is arbitrary, we are done.\hfill$\blacksquare$
\vskip 0.1in
\noindent
{\it Proof of Theorem 4.1}. Identify $A=L^{\infty}(X)$ and assume by contradiction that $A$ and $B$ are not unitarily conjugate. By Proposition 4.5 we can find  $X_0\subset X$ of positive measure, equivalence relations $\Cal S\subset\Cal T\subset \Cal R_{|X_0}$ and a measurable partition $\{X_k\}_{k\geqslant 1}$ of $X_0$  such that  
(1) $\Cal S$ is hyperfinite and has infinite orbits,
(2) $\Cal S$ is q--normal in $\Cal T$, and 
(3) almost every $\Cal R_{|X_k}$--class contains only finitely many $\Cal T_{|X_k}$--classes, for all $k\geqslant 1$.

It is easy to see that (3) implies that $\Cal T$ is q--normal in $\Cal R_{|X_0}$. Since $\Cal S$ is q--normal in $\Cal T$, by applying Lemma 4.6 twice we get that $\Cal C(\Cal R_{|X_0})\leqslant \Cal C(\Cal S)$. This is a contradiction because the induction formula [Ga99, Proposition II.6.] gives that $\Cal C(\Cal R_{|X_0})=1+\mu(X_0)^{-1}(\Cal C(\Cal R)-1)>1$, while the fact that $\Cal S$ is hyperfinite implies that $\Cal C(\Cal S)\leqslant 1$ (see [Ga99, Proposition III.3.]).\hfill$\blacksquare$

\vskip 0.1in
\noindent
{\it Remark.} Consider the usual action SL$_2(\Bbb Z)\curvearrowright (\Bbb T^2,\lambda^2)$ and let $M=L^{\infty}(\Bbb T^2)\rtimes$ SL$_2(\Bbb Z)$. Then by using the results of the last two sections and [Oz08] we can already show that $M$ has a unique group measure space Cartan subalgebra. Indeed, assume that $M=L^{\infty}(Y)\rtimes\Lambda$, for some free ergodic p.m.p. action $\Lambda\curvearrowright (Y,\nu)$. Firstly, by Theorem 3.1 we get that $L^{\infty}(X)\prec_{M}L^{\infty}(Y)\rtimes\Sigma$, for a subgroup $\Sigma<\Lambda$ which is either amenable or of the form $\Sigma=\cup_{n\geqslant 1}C(\Lambda_n)$, for a decreasing family $\{\Lambda_n\}_{n\geqslant 1}$ of infinite subgroups of $\Lambda$. Secondly, since $M$ is solid [Oz08], we deduce that $\Sigma$ must be amenable in either case. Finally, by  Theorem 4.2 we conclude that $L^{\infty}(X)$ and $L^{\infty}(Y)$ are unitarily conjugate.

\vskip 0.2in
\head \S 5. {Proof of Theorem 1}.\endhead

\vskip 0.2in
In this section we combine the results of the previous section to prove Theorem 1 and more generally:

\proclaim {Theorem 5.1} Let $\Gamma$ be an infinite countable group with  $\beta_1^{(2)}(\Gamma)>0$. Let  $\Gamma\curvearrowright (X,\mu)$ be a free ergodic rigid p.m.p. action. Let $s>0$ and
denote $M=L^{\infty}(X)\rtimes\Gamma$. 
\vskip 0.03in
\noindent
If $\Lambda\curvearrowright (Y,\nu)$ is any free ergodic p.m.p. action  such that 
 $M^s=L^{\infty}(Y)\rtimes\Lambda$,  then we can find a unitary $u\in M^s$ such that 
$uL^{\infty}(X)^su^*=L^{\infty}(Y)$. 
\endproclaim
\noindent
{\it Proof.} Consider a group measure space decomposition $M^{s}=B\rtimes\Lambda$, for  $s>0$. Let $n\geqslant s$ be an integer and $p\in D_n(\Bbb C)\otimes L^{\infty}(X)$ be a projection of trace $\frac{s}{n}$. Identify $M^s=p(\Bbb M_n(\Bbb C)\otimes M)p$ and $L^{\infty}(X)^s=p(D_n(\Bbb C)\otimes L^{\infty}(X))p$.  Let $\frac{\Bbb Z}{n\Bbb Z}$ act on itself by addition and endow
 $\tilde X=X\times\frac{\Bbb Z}{n\Bbb Z}$  with the diagonal action of $\tilde\Gamma=\Gamma\times\frac{\Bbb Z}{n\Bbb Z}.$ 
Then $\beta_1^{(2)}(\tilde\Gamma)>0$, the action $\tilde\Gamma\curvearrowright\tilde X$ is free ergodic rigid p.m.p. and we have that $\Bbb M_n(\Bbb C)\otimes M=L^{\infty}(\tilde X)\rtimes\tilde\Gamma$ and $D_n(\Bbb C)\otimes L^{\infty}(X)=L^{\infty}(\tilde X)$. 
Thus, after replacing $\Gamma, X$ with $\tilde\Gamma, \tilde X$, we may assume that $s\leqslant 1$, i.e. $pMp=B\rtimes\Lambda$, for a projection $p\in L^{\infty}(X)$. 

\vskip 0.05in
Since the action $\Gamma\curvearrowright X$ is rigid, the inclusion $L^{\infty}(X)p\subset pMp$ has the relative property (T) ([Po01, Proposition 4.7]). Also, since $\Gamma$ has positive first $\ell^2$--Betti number, it admits an unbounded cocycle $b:\Gamma\rightarrow\ell^2_{\Bbb R}\Gamma$ ([PT07, Corollary 2.4]).  Altogether, by applying Theorem 3.1 we are in one of the following two situations:

\vskip 0.02in
\noindent
{\bf Case 1.} $L^{\infty}(X)p\prec_{pMp}B\rtimes\Lambda_0$, for an amenable subgroup $\Lambda_0$ of $\Lambda$.
\vskip 0.02in
\noindent
{\bf Case 2.} $L^{\infty}(X)p\prec_{pMp}B\rtimes (\cup_{n\geqslant 1}C(\Lambda_n))$, for a decreasing sequence $\{\Lambda_n\}_{n\geqslant 1}$ of non--amenable subgroups of $\Lambda$.

\vskip 0.05in
In the first case, Theorem 4.2 gives the conclusion. Thus, we may assume that we are in the second case.
If the group $\cup_{n\geqslant 1}C(\Lambda_n)$ is amenable, then we are again in the first case. So, we may additionally assume that $\cup_{n\geqslant 1}C(\Lambda_n)$ is non--amenable. It follows that $C(\Lambda_n)$ is  non--amenable, for some $n\geqslant 1$.

Let $\tilde M\supset M$ and the automorphisms $\{\alpha_t\}_{t\in\Bbb R}$ of $\tilde M$  be as defined in Section 2. Since $C(\Lambda_n)$ is non--amenable, $L(C(\Lambda_n))$ has no amenable direct summand and Lemma 2.2 implies that $\alpha_t\rightarrow id$ uniformly on $(L\Lambda_n)_1$.
Since $\Lambda_n$ is non--amenable, [Po03, Theorem 2.1 and Corollary 2.3] provides a sequence $g_k\in\Lambda_n$ such that $||E_{L^{\infty}(X)}(xv_{g_k}y)||_2\rightarrow 0$, for all $x,y\in M$ (here $\{v_g\}_{g\in\Lambda}\in B\rtimes\Lambda$ denote the canonical unitaries).

Further, applying Theorem 2.4 to $\{v_{g_k}\}_{k\geqslant 1}$ gives that $\alpha_t\rightarrow id$ uniformly on $(B)_1$. Finally, Theorem 2.3 implies that either $B\prec_{M}L^{\infty}(X)$ or $\alpha_t\rightarrow id$ uniformly on $(pMp)_1$. In the first case Lemma 1.3 yields that $B$ and $L^{\infty}(X)p$ are unitarily conjugate while in the second case, Lemma 2.1  implies that $b$ is bounded, a contradiction.
\hfill$\blacksquare$

\vskip 0.05in
\noindent
{\it Remark.}  Let us recall Ozawa and Popa's examples of $\Cal H\Cal T$ factors with two non--conjugate Cartan subalgebras ([OP08]) and explain why Theorem 5.1 does not apply to them.
Let $p_1,p_2,...$ be prime numbers and define $G=\cup_{n\geqslant 1}\{z\in\Bbb T|z^{p_1p_2\cdots p_n}=1\}$. Then $G^2<\Bbb T^2$ is an SL$_2(\Bbb Z)$--invariant subgroup and $\Gamma=G^2\rtimes$ SL$_2(\Bbb Z)$ has Haagerup's property. 
Also, the action $\Gamma\curvearrowright (\Bbb T^2,\lambda^2)$ (where $G^2$ and SL$_2(\Bbb Z)$ act on $\Bbb T^2$ by translations and automorphisms, respectively) is free ergodic and rigid. Thus, $M=L^{\infty}(\Bbb T^2)\rtimes\Gamma$ is an $\Cal H\Cal T$ factor.
Moreover, as shown in [OP08] and [PV09, Section 5.5], $L(G^2)$ is a group measure space Cartan subalgebra of $M$ which is not conjugate to $L^{\infty}(\Bbb T^2)$. 

Since $\Gamma$ has an infinite normal abelian subgroup,  [CG86] gives that $\beta_1^{(2)}(\Gamma)=0$, showing why Theorem 5.1 does not apply to $M$.  
\vskip 0.05in

\vskip 0.02in
\head\S 6. A strong rigidity result and applications.\endhead
\vskip 0.2in
Let $\Gamma$ be a countable group with positive first $\ell^2$--Betti number. Then a far--reaching conjecture of Chifan, Peterson, Popa and the author predicts that any II$_1$ factor $L^{\infty}(X)\rtimes\Gamma$,  arising from a free ergodic p.m.p. action $\Gamma\curvearrowright (X,\mu)$,  has a unique  Cartan subalgebra (see [Po09]). Chifan and Peterson proved that  if $\Gamma$ admits a non--amenable subgroup with the relative property (T), then $L^{\infty}(X)\rtimes\Gamma$ has a unique group measure space Cartan subalgebra ([CP10, Theorem 7.4]).

In this section,  we  weaken the rigidity assumption on $\Gamma$ by requiring that $\Gamma$ does not have Haagerup's property and show that a lot can still be said about the group measure space decompositions of $L^{\infty}(X)\rtimes\Gamma$. 
Although, in general,  we cannot conclude that $L^{\infty}(X)\rtimes\Gamma$ has a unique group measure Cartan subalgebra, we deduce that this is the case if $\Gamma\curvearrowright (X,\mu)$ is a  {\it solid} action (see Corollary 6.4).
\proclaim {Theorem 6.1} Let $\Gamma\curvearrowright (X,\mu)$ be a free ergodic p.m.p. action and
denote $M=L^{\infty}(X)\rtimes\Gamma$. Assume that 
 $\beta_1^{(2)}(\Gamma)>0$ and $\Gamma$ does not have Haagerup's property.
Let $\Lambda\curvearrowright (Y,\nu)$ be a free ergodic p.m.p. action  such that 
 $M^s=L^{\infty}(Y)\rtimes\Lambda$, for some $s>0$. Suppose that $L^{\infty}(X)^s$ and $L^{\infty}(Y)$ are not unitarily conjugate. Then we have that

\vskip 0.03in
\noindent
(1) $\Lambda$ does not have Haagerup's property.

\vskip 0.02in
\noindent
(2) We can find an infinite abelian subgroup $\Delta_0<\Lambda$ such that $L\Delta_0\prec_{M^s}L^{\infty}(X)^s$ and the centralizer of $\Delta_0$ in $\Lambda$ is non--amenable. 

\vskip 0.02in
\noindent
(3) For every $h\in\Lambda$, we can find a finite index subgroup $\Delta_1<\Delta_0$ such that the groups  $h\Delta_1 h^{-1}$ and $\Delta_1$ commute.

\vskip 0.02in\noindent
(4)  $\beta_1^{(2)}(\Lambda)=0$.
\endproclaim

\noindent
{\it Remark.} If $L^{\infty}(X)^s$ and $L^{\infty}(Y)$ {\it are} unitarily conjugate, then the involved actions are stably orbit equivalent. Since Haagerup's property is invariant under stable orbit equivalence (see e.g. [Po01, Corollary 2.5 and Proposition 3.1]),  we also get that $\Lambda$ does not have Haagerup's property.
\vskip 0.05in

In the proof of Theorem 6.1 we will need the following lemma due to Houdayer, Popa and Vaes.

\proclaim  {Lemma 6.2 [HPV10]} Let $(A,\tau)$ be a tracial von Neumann algebra and $\Gamma\curvearrowright (A,\tau)$ be a trace preserving action. Denote $M=A\rtimes\Gamma$ and let $B\subset pMp$ be a regular von Neumann subalgebra. Assume that $B\prec_{M}A\rtimes\Sigma$, for some subgroup $\Sigma$ of $\Gamma$. 

\noindent
Denote by $\Delta$ the subgroup of $\Gamma$ generated by all $g\in\Gamma$ such that $g\Sigma g^{-1}\cap\Sigma$ is infinite.

\vskip 0.03in
\noindent
If $B\nprec_{M}A$, then $\Delta$ has finite index in $\Gamma$.
\endproclaim

\noindent
{\it Proof.} By Section 4 in [HPV10], given a subgroup $\Sigma<\Gamma$, we can find a projection $z(\Sigma)\in M$
such that $z(\Sigma)\not=0$ iff $B\prec_{M}A\rtimes\Sigma$ and $z(g\Sigma g^{-1})=u_gz(\Sigma)u_g^*$, for all $g\in\Gamma$. Moreover, by  [HPV10, Proposition 6], $z(\Sigma\cap\Sigma')=z(\Sigma)z(\Sigma')$, for any subgroup $\Sigma'<\Gamma$.

Assume by contradiction that $\Delta$  has infinite index in $\Gamma$. Then we can find $\{g_i\}_{i\geqslant 1}\subset\Gamma$ such that $g_i\Sigma g_i^{-1}\cap g_j\Sigma g_j^{-1}$ is finite, for every $i,j\geqslant 1$. Since $B\nprec_{M}A$, it follows that $z(g_i\Sigma g_i^{-1}\cap g_j\Sigma g_j^{-1})=0$, for every $i,j\geqslant 1$. By using the above formulas we derive that the projections $\{u_{g_i}z(\Sigma)u_{g_i}^*\}_{i\geqslant 1}$ are mutually orthogonal. Since $z(\Sigma)\not=0$, this leads to a contradiction.\hfill$\blacksquare$ 

\vskip 0.1in

\noindent
{\it Proof of Theorem 6.1.} 
 By reasoning as in the beginning of Section 5, we can reduce to the case $s\leqslant 1$. Therefore, we may assume that $pMp=B\rtimes\Lambda$, where $p\in A=L^{\infty}(X)$ is a projection and $B=L^{\infty}(Y)$. Denote by $\{u_g\}_{g\in\Gamma}\subset M$ and $\{v_h\}_{h\in\Lambda}\subset pMp$ the canonical unitaries.
Since $Ap$ and $B$ are not unitarily conjugate and $\beta_{1}^{(2)}(\Gamma)>0$, Theorem 4.2 implies the following fact that we will use repeatedly:
\vskip 0.05in
\noindent
 {\bf Fact.} If $A\prec_{M}B\rtimes\Sigma$, for a subgroup $\Sigma<\Lambda$, then  $\Sigma$ is non--amenable. 

\noindent
Similarly, if $B\prec_{M}A\rtimes\Sigma$,   for a subgroup $\Sigma<\Gamma$, then  $\Sigma$ is non--amenable. 
\vskip 0.05in
The  proof of Theorem 6.1 is split between five claims, all of which, with the exception of Claim 2, prove one of the conditions (1)--(4) from the conclusion. 

\vskip 0.05in
\noindent
{\bf Claim 1.} $\Lambda$ does not have Haagerup's property.
\vskip 0.05in
\noindent
{\it Proof of Claim 1.} Assuming by contradiction that $\Lambda$ has Haagerup's property, 
 we can find a sequence $\phi_n:\Lambda\rightarrow \Bbb C$ of positive definite functions such that $\phi_n(h)\rightarrow 1$, for all $h\in\Lambda$, and $\phi_n\in c_0(\Lambda)$, for all $n\geqslant 1$. As $M$ is a factor there are partial isometries $w_1,..,w_k\in M$ such that $w_iw_i^*\leqslant p$, for all $i,$ and $\sum_{i=1}^kw_i^*w_i=1$.
For $n\geqslant 1$, we define
\vskip 0.02in
$\bullet$ $\Phi_n:pMp\rightarrow pMp$ by \hskip 0.05in $\Phi_n(x)=\sum_{h\in\Lambda}\phi_n(h)b_hv_h$, for all $x=\sum_{h\in\Lambda}b_hv_h\in pMp$,	

$\bullet$  $\Psi_n:M\rightarrow M$ by letting \hskip 0.05in $\Psi_n(x)=\sum_{i,j=1}^kw_i^*\Phi_n(w_ixw_j^*)w_j$, for all $x\in M$, and

$\bullet$ $\psi_n:\Gamma\rightarrow\Bbb C$ by letting \hskip 0.05in $\psi_n(g)=\tau(\Psi_n(u_g)u_g^*)$, for all $g\in\Gamma$.
\vskip 0.05in
Then $\psi_n$ are positive definite functions and $\psi_n(g)\rightarrow 1$, for al $g\in\Gamma$. Since $\Gamma$ does not have Haagerup's property, 
[Pe09, Lemma 2.6] provides $n_0\geqslant 1$ and an infinite sequence $\{g_m\}_{m\geqslant 1}\subset\Gamma$ such that $\inf_{m}|\psi_{n_0}(g_m)|\geqslant\frac{1}{2}$. Thus, we have  $\inf_{m}||\Psi_{n_0}(u_{g_m})||_2\geqslant\frac{1}{2}$.
On the other hand, it is easy to see that $\Psi_{n_0}$ is ``compact over $B$": if a sequence $x_m\in (M)_1$ satisfies $||E_B(yx_mz)||\rightarrow 0$, for all $y,z\in M$, then $||\Psi_{n_0}(x_m)||_2\rightarrow 0$.

 The last two facts imply that, after replacing $\{g_m\}_{m\geqslant 1}$ with a subsequence, we can find $y,z\in M$ such that $\inf_{m}||E_B(yu_{g_m}z)||_2>0$. Moreover, we may clearly assume  that $y,z\in (A)_1$. 
For $m\geqslant 1$, let  $b_m=E_B(yu_{g_m}z)$.
 Since  $b_m\in B$ and $a_m:=(u_{g_m}z^*u_{g_m}^*)y\in (A)_1$, we get that $||b_m||_2^2=\tau(b_mz^*u_{g_m}^*y^*)=\tau(b_mu_{g_m}^*a_m)\leqslant ||E_A(b_mu_{g_m}^*)||_2$. 
Since $\inf_{m}||b_m||_2>0$, it follows that $\inf_{m}||E_A(b_mu_{g_m}^*)||_2>0$.

By applying Lemma 3.4 we get that $B\prec_{M}A\rtimes\Sigma$, where $\Sigma=\cup_{m\geqslant 1}C(\Gamma_m)$, for some decreasing sequence $\{\Gamma_m\}_{m\geqslant 1}$  of infinite subgroups of $\Gamma$. 

To reach a contradiction it suffices to show that any cocycle $c:\Gamma\rightarrow\ell^2\Gamma$ for the regular representation $\pi:\Gamma\rightarrow\ell^2\Gamma$ is inner.  
Since $\Sigma$ is non--amenable (by the above Fact),  $C(\Gamma_{m_0})$ is non--amenable for some $m_0\geqslant 1$. 
By Lemma 2.5 (1) we can find $\xi\in\ell^2\Gamma$ such that $c(g)=\pi(g)\xi-\xi,$ for all $g\in\Gamma_{m_0}$. Let $\Gamma_0\subset\Gamma$ be the subgroup of all $g\in\Gamma$ such that $c(g)=\pi(g)\xi-\xi$. 
If $m\geqslant m_0$, then $\Gamma_m\subset\Gamma_{m_0}\subset\Gamma_0$. Since $\Gamma_m$ is infinite by Lemma 2.5 (2) it follows that $C(\Gamma_m)\subset \Gamma_0$ and thus $\Sigma\subset \Gamma_0$. 

Now, denote by $\Delta$ the subgroup $\Gamma$ generated by all $g\in\Gamma$ for which $g\Sigma g^{-1}\cap\Sigma$ is infinite. Note that if $g\Sigma g^{-1}\cap\Sigma$ is infinite, then $g\Gamma_0g^{-1}\cap \Gamma_0$ is infinite and therefore $g\in\Gamma_0$ (by Lemma 2.5 (2)). This shows that $\Delta\subset\Gamma_0$.
On the other hand, since $B\prec_{M}A\rtimes\Sigma$ but $B\nprec_{M}A$, Lemma 6.2 implies that $\Delta$ has finite index in $\Gamma$. 
Thus, $\Gamma_0$ has finite index in $\Gamma$ and by applying Lemma 2.5 (2) again we conclude that $\Gamma_0=\Gamma$. In other words, $c$ is inner, as claimed.\hfill$\square$
\vskip 0.05in

Next, let $b:\Gamma\rightarrow\ell^2_{\Bbb R}\Gamma$ be an unbounded cocycle for the left regular representation. Let $\tilde M\subset M$ and $\{\alpha_t\}_{t\in\Bbb R}$ be defined as in Section 2.
By using Claim 1 we deduce: 
\vskip 0.05in
\noindent
{\bf Claim 2.}  There exist an infinite sequence $\{h_n\}_{n\geqslant 1}\subset \Lambda$ and $x\in M$ such that $\inf_{n}||E_A(xv_{h_{n}})||_2>0$.
\vskip 0.05in
\noindent
{\it Proof of Claim 2.}
For $t\in\Bbb R$, define a positive definite function $\phi_t:\Lambda\rightarrow\Bbb C$ through the formula $\phi_t(h)=\tau(\alpha_t(v_h)v_h^*)$, for $h\in\Lambda$.
Then $\phi_t(h)\nearrow \tau(p)$, as $t\rightarrow 0$, for all $h\in\Lambda$. Since $\Lambda$ does not have Haagerup's property, by [Pe09, Lemma 2.6] we can find  an infinite sequence $\{h_n\}_{n\geqslant 1}\subset \Lambda$  such that $\sup_{n\geqslant 1}|\tau(p)-\phi_t(h_n)|\rightarrow 0$, as $t\rightarrow 0$. It follows that $\alpha_t\rightarrow id$
uniformly on $\{v_{h_n}\}_{n\geqslant 1}$. 

If the claim is false, then
$||E_A(xv_{h_n})||_2\rightarrow 0$, for all $x\in M$. Thus, $||E_A(xv_{h_n}y)||_2\rightarrow 0$, for all $x,y\in M$. Since $\{v_{h_n}\}_{n\geqslant 1}$ normalize $B$, Theorem 2.4 implies that $\alpha_t\rightarrow id$ uniformly on $(B)_1$. Since $B\nprec_{M}A$, Theorem 2.3 gives that $\alpha_t\rightarrow id$ uniformly on $(pMp)_1$. But then Lemma 2.1 would imply that $b$ is bounded, a contradiction.\hfill$\square$
\vskip 0.02in
 Let $\{h_n\}_{n\geqslant 1}$ and $x\in M$ as given by Claim 2. Since $E_A(xv_{h_n})=E_A(pxpv_{h_n})$, we may assume that $x\in pMp=B\rtimes\Lambda$. By replacing $h_n$ with a subsequence we can assume that $x=bv_h$, for some $b\in (B)_1$ and $h\in\Lambda$. Finally, by replacing $h_n$ with $hh_n$, we can assume that $\inf_{n}||E_A(bv_{h_n})||_2>0$, for some $b\in (B)_1$.

\vskip 0.05in
\noindent
{\bf Claim 3.} There exists an  infinite abelian subgroup $\Delta_0<\Lambda$ with non--amenable centralizer such that $(L\Delta_0)q\prec_{M}A$, for every non--zero projection $q\in L\Delta_0'\cap B$.
\vskip 0.05in
\noindent
{\it Proof of Claim 3.} For every $n\geqslant 1$, denote $a_n=E_A(bv_{h_n})$. Then $a_n\in (Ap)_1$ and $\inf_n||a_n||_2>0$.  Also, since $a_n\in A$ and $b\in (B)_1$, we get that $$||a_n||_2^2=\tau(a_nv_{h_n}^*b^*)\leqslant ||E_B(a_nv_{h_n}^*)||_2.$$ By combining the last two inequalities we derive that $\inf_{n}||E_B(a_nv_{h_n}^*)||_2>0$. Since $a_n\in (Ap)_1$ and $h_n\rightarrow\infty$, Lemma 3.4 implies that $Ap\prec_{M}B\rtimes\Sigma$, where $\Sigma=\cup_{m\geqslant 1}C(\Lambda_m)$, for some decreasing sequence $\{\Lambda_m\}_{m\geqslant 1}$ of infinite subgroups of $\Lambda$.

Next, by the above Fact, $\Sigma$ is non--amenable. Thus, $C(\Lambda_{m_0})$ is non--amenable for some $m_0\geqslant 1$.  Put $\Delta=\Lambda_{m_0}$. Lemma 2.2 then gives that $\alpha_t\rightarrow id$ uniformly on $(L\Delta)_1$.
We claim that $(L\Delta) q\prec_{M}A$, for every non--zero projection $q\in (L\Delta)'\cap B$.

 Otherwise, by
[Po03, Theorem 2.1 and Corollary 2.3] we can find a sequence $\lambda_i\in\Delta$ such that $||E_A(xv_{\lambda_i}qy)||\rightarrow 0$, for all $x,y\in M$. Note that $v_{\lambda_i}q\in \Cal U(qMq)$ normalizes $Bq$, for all $i\geqslant 1$, and that $\alpha_t\rightarrow id$ uniformly on $\{v_{\lambda_i}q\}_{i\geqslant 1}$. But then Theorem 2.4 would give that $Bq\prec_{M}A$, a contradiction.

Since $L\Delta\prec_{M}A$, we get that $\Delta$ is virtually abelian. Let $\Delta_0<\Delta$ be a finite index abelian subgroup. Since $\alpha_t\rightarrow id$ uniformly on $(L\Delta_0)_1$, arguing as in the previous paragraph shows that $(L\Delta_0) q\prec_{M}A$, for every non--zero projection $q\in (L\Delta_0)'\cap B$.\hfill$\square$

\vskip 0.05in
\noindent
{\bf Claim 4.} For every $h\in\Lambda$, we can find a finite index subgroup $\Delta_1<\Delta_0$ such that the groups  $h\Delta_1 h^{-1}$ and $\Delta_1$ commute.
\vskip 0.05in
\noindent
{\it Proof of Claim 4.}
 Let $\Omega_0$ be the group of $k\in\Lambda$ for which the set $\{\lambda k\lambda^{-1}|\lambda\in\Delta_0\}$ is finite, i.e. such that $k$ commutes with a finite index subgroup of $\Delta_0$. Then $\Delta_0\subset\Omega_0$ and $(L\Delta_0)'\cap B\rtimes\Lambda\subset B\rtimes\Omega_0$. 

Now, let $r\in (B\rtimes\Omega_0)'\cap pMp$ be a non--zero projection.  Since $\Delta_0\subset\Omega_0$ and $B\subset pMp$ is maximal abelian, it follows that $r\in (L\Delta_0)'\cap B$. By Claim 3 we get that $(L\Delta_0) r\prec_{M}A$. Since $A\subset M$ is a Cartan subalgebra, it follows that $(L\Delta_0) r\prec_{pMp}Ap$. By taking relative commutants we get that $Ap\prec_{pMp}(B\rtimes\Omega_0)r$ ([Va07, Lemma 3.5]).

 Since $Ap\subset pMp=B\rtimes\Lambda$ is regular, [HPV10, Corollary 7] implies that $Ap\prec_{pMp}B\rtimes(h\Omega_0 h^{-1}\cap\Omega_0),$ for every $h\in\Lambda$. Fix $h\in\Lambda$. Then the Fact from the beginning of the proof gives that $h\Omega_0 h^{-1}\cap\Omega_0$ is non--amenable. Let  $\Omega<\Omega_0$ be a finitely generated subgroup such that  $\Sigma:=h\Omega h^{-1}\cap\Omega$ is  also non--amenable. Since every element of $\Omega_0$ commutes with a finite index subgroup of $\Delta_0$ and $\Omega$ is finitely generated, we can find a finite index subgroup $\Delta<\Delta_0$ which commutes with $\Omega$. 

Let $\Upsilon$ be the subgroup of $\Lambda$ generated by $h\Delta h^{-1}$ and  $\Delta$.
Then $\Sigma$ and $\Upsilon$ commute.
 Since $\Sigma$ is non--amenable, arguing as in the proof of Claim 3 gives that $\Upsilon$ is virtually abelian. The claim now follows easily.\hfill$\square$
\vskip 0.05in
\noindent
{\bf Claim 5.} $\beta_1^{(2)}(\Lambda)=0$.  
\vskip 0.05in
\noindent
{\it Proof of Claim 5.}
Let $c:\Lambda\rightarrow\ell^2\Lambda$ be a cocycle for the regular representation. Since by Claim 3, $\Delta_0$ has non--amenable centralizer in $\Lambda$, Lemma 2.5 (1) provides a vector $\xi\in\ell^2\Lambda$  such that $c(g)=\pi(g)\xi-\xi$, for all $g\in\Delta_0$. 

Let $\Lambda_0<\Lambda$ the subgroup of $g\in\Lambda$ such that $c(g)=\pi(g)\xi-\xi$.
Let $h\in\Lambda$. By Claim 4 there is finite index subgroup $\Delta_1<\Delta_0$ such that  $h^{-1}\Delta_1h$ and $\Delta_1$ commute. Since $\Delta_1$ is infinite and $\Delta_1<\Lambda_0$, Lemma 2.5 (2) gives that $h^{-1}\Delta_1h<\Lambda_0$. Thus $\Delta_1<h\Lambda_0 h^{-1}\cap \Lambda_0$ and Lemma 2.5 (2) yields that $h\in\Lambda_0$. This shows that $\Lambda_0=\Lambda$, i.e. $c$ is inner.  This finishes the proofs of the claim and of the theorem. \hfill$\blacksquare$ 
\vskip 0.05in
We can now deduce corollaries 4 and 5 stated in the introduction.

\proclaim {Corollary 6.3} Let $\Gamma$ be a countable group such that $\beta_1^{(2)}(\Gamma)\in (0,+\infty)$ and $\Gamma$ does not have Haagerup's property.
 Let $\Gamma\curvearrowright (X,\mu)$ be any free ergodic p.m.p. action.
 
\vskip 0.03in
\noindent
Then the II$_1$ factor $M=L^{\infty}(X)\rtimes\Gamma$ has trivial fundamental group, $\Cal F(M)=\{1\}$.
\endproclaim

Note that under the stronger assumption that $\Gamma$ has a non--amenable subgroup with the relative property (T)  this result also follows from [Va10b, Theorem 1.3].

\vskip 0.05in
\noindent
{\it Proof.} For $t\in\Cal F(M)$, let $\theta:M^t\rightarrow M$ be an isomorphism. Then we can find a unitary $u\in M$ such that $u\theta(L^{\infty}(X)^t)u^*=L^{\infty}(X)$. Indeed, otherwise by Theorem 6.1 we would get that $\beta_1^{(2)}(\Gamma)=0,$ a contradiction.
Thus, if $\Cal R$ denotes the equivalence relation induced by the action $\Gamma\curvearrowright (X,\mu)$, then $\Cal R^t\cong\Cal R$. This shows that $\Cal F(M)=\Cal F(\Cal R)$. 

On the other hand, [Ga01, Corollaire 3.17] gives that $\beta_{1}^{(2)}(\Cal R)=\beta_{1}^{(2)}(\Gamma)\in (0,+\infty)$. By applying [Ga01, Corollaire 5.7] we deduce that $\Cal F(\Cal R)=\{1\}$, thus $\Cal F(M)=\{1\}$.\hfill$\blacksquare$

\proclaim {Corollary 6.4}  Let $\Gamma$ be a countable group such that $\beta_1^{(2)}(\Gamma)>0$ and $\Gamma$ does not have Haagerup's property. Assume that one of the following two conditions holds true:

\noindent (1) $\Gamma\curvearrowright (X,\mu)=(X_0^I,\mu_0^I)$ is a free, generalized Bernoulli action, where  $(X_0,\mu_0)$ is a non--trivial probability space and $\Gamma\curvearrowright I$ is an action  with amenable stabilizers.

\noindent (2) $\Gamma\curvearrowright (X,\mu)$ is a free ergodic p.m.p {\it solid} action, i.e. the relative commutant $Q'\cap L^{\infty}(X)\rtimes\Gamma$ is amenable,  for any diffuse von Neumann subalgebra $Q\subset L^{\infty}(X)$.
\vskip 0.03in
\noindent
If $\Lambda\curvearrowright (Y,\nu)$ is any free ergodic p.m.p. action such that $M^t=L^{\infty}(Y)\rtimes\Lambda$, for some $t>0$, then we can find a unitary element $u\in M^t$ such that $uL^{\infty}(X)^tu^*=L^{\infty}(Y)$. 
\endproclaim
\noindent
{\it Proof.} Firstly, [CI08, Theorem 7] gives that (1) $\Longrightarrow$ (2), so we can assume that (2) is satisfied.
 Now, suppose by contradiction that the conclusion is false. Then by Theorem 6.1 we can find an infinite subgroup $\Delta_0<\Lambda$ such that its centralizer is non--amenable and $L\Delta_0\prec_{M^t}L^{\infty}(X)^t$. It follows that we can find a diffuse von Neumann subalgebra $D\subset L^{\infty}(X)^t$ such that $D'\cap M^t$ is non--amenable. This however contradicts the assumption that $\Gamma\curvearrowright (X,\mu)$ is solid. 
\hfill$\square$

\head References\endhead
\item {[BO08]} N.P. Brown, N. Ozawa: {\it C$^*$--algebras and finite-dimensional approximations,} Graduate
Studies in Mathematics, 88. American Mathematical Society, Providence, RI,
2008. xvi+509 pp.
\item {[BV97]} M. Bekka, A. Valette: {\it Group cohomology, harmonic functions and the first $L^2$--Betti number}, Potential Anal. {\bf 6} (1997), no. 4, 313�-326.
\item {[Bu91]} M. Burger: {\it Kazhdan constants for} SL$(3,\Bbb Z)$, J. Reine Angew. Math. {\bf 413} (1991), 36--67.
\item {[CFW81]} A. Connes, J. Feldman, B. Weiss: {\it An amenable equivalence relation is generated
by a single transformation},  Ergodic. Th. and Dynam. Sys {\bf 1} (1981), no. 4, 431--450.
\item {[CG86]} J. Cheeger, M. Gromov: {\it $L_2$--cohomology and group cohomology}, Topology {\bf 25} (1986), no. 2, 189�-215.
\item {[CI08]} I. Chifan, A. Ioana: {\it Ergodic subequivalence relations induced by a Bernoulli action}, Geom. Funct. Anal. Vol. 20 (2010), 53--67. 
\item {[CP10]} I. Chifan, J. Peterson: {\it Some unique group-measure space decomposition results}, preprint arXiv:1010.5194.
\item {[CS11]} I. Chifan, T. Sinclair: {\it On the structural theory of II$_1$ factors of negatively curved groups}, preprint arXiv:1103.4299.

\item {[FM77]} J. Feldman, C.C. Moore: {\it Ergodic equivalence relations, cohomology, and von Neumann algebras, II}, Trans. Amer. Math. Soc. {\bf 234} (1977), 325--359.
\item {[FV10]} P. Fima, S. Vaes: {\it HNN extensions and unique group measure space decomposition of II$_1$ factors}, Trans. Amer. Math. Soci. {\bf 364} (2012), 2601--2617
\item {[Fu09]} A. Furman: {\it A survey of Measured Group Theory}, Geometry, Rigidity, and Group Actions, 296--374, The University of Chicago Press, Chicago and London, 2011.
\item {[Ga99]} D. Gaboriau: {\it Co$\hat{u}$t des relations d'equivalence et des groupes. (French) [Cost of
equivalence relations and of groups]} Invent. Math. {\bf 139} (2000), no. 1, 41�-98.
\item {[Ga01]} D. Gaboriau: {\it Invariants L$^2$ de relations d'equivalence et de groupes}, Publ.
Math. Inst. Hautes \'Etudes Sci., {\bf 95} (2002), 93--150.
\item {[Ga08]} D. Gaboriau: {\it Relative Property (T) Actions and Trivial Outer Automorphism }

{\it Groups}, J. Funct. Anal. {\bf 260} (2011), no. 2, 414--427.
\item {[Ga10]} D. Gaboriau: {\it Orbit Equivalence and Measured Group Theory},   Proceedings of the ICM (Hyderabad, India, 2010), Vol. III, Hindustan Book Agency (2010), 1501--1527. 
\item {[HPV10]} C. Houdayer, S. Popa, S. Vaes: {\it A class of groups for which every action is W$^*$--superrigid}, preprint arXiv:1010.5077, to appear in Groups Geom. Dyn.,
\item {[IKT08]} A. Ioana, A. S. Kechris, T. Tsankov: {\it Subequivalence relations and positive-definite functions},   Groups Geom. Dyn.,  Volume 3, Issue 4, (2009), 579--625.
\item {[Io07]} A. Ioana: {\it Orbit inequivalent actions for groups containing a copy of $\Bbb F_2$},   Invent. Math. {\bf 185} (2011), 55--73.
\item {[Io09]} A. Ioana: {\it Relative property (T) for the subequivalence relations induced by the action of} SL$(2,\Bbb Z)$ {\it on} $\Bbb T^2$,
Advances in Math.  {\bf 224} (2010), 1589--1617. 
 \item {[Io10]} A. Ioana: {\it W$^*$--superrigidity for Bernoulli actions of property (T) groups}, J. Amer. Math. Soc. {\bf 24} (2011), 1175--1226.
\item {[IPV10]} A. Ioana, S. Popa, S. Vaes: {\it A class of superrigid group von Neumann algebras}, preprint arXiv:1007.1412.
 \item {[IS10]} A. Ioana, Y. Shalom: {\it Rigidity for equivalence relations on homogeneous spaces}, preprint  arXiv:1010.3778, to appear in Groups Geom. Dyn.
\item {[Ku51]} M. Kuranishi, {\it On everywhere dense embedding of free groups in Lie groups}, Nagoya
Math. J. {\bf 2} (1951), 63-�71.
\item {[MS06]} N. Monod, Y. Shalom: {\it Orbit equivalence rigidity and bounded cohomology}, Ann. of Math. (2), {\bf 164} (2006), no. 3, 825�-878.
\item {[MvN36]} F. Murray, J. von Neumann: {\it On rings of operators,} Ann. of Math. {\bf 37} (1936), 116--229.
\item {[Oz08]} N. Ozawa: {\it An example of a solid von Neumann algebra}, Hokkaido Math. J., {\bf 38} (2009), 557--561.
\item {[OP07]} N. Ozawa, S. Popa: {\it On a class of II$_1$ factors with at most one Cartan subalgebra}, Ann. of Math. (2), {\bf 172} (2010), 713--749.
\item {[OP08]} N. Ozawa, S. Popa: {\it On a class of II$_1$ factors with at most one Cartan subalgebra, II}, Amer. J. Math., {\bf 132} (2010), 841--866. 
\item {[Pe06]} J. Peterson: {\it $L^2$--rigidity in von Neumann algebras},  Invent. Math.  {\bf 175}  (2009),  no. 2, 417--433.
\item {[Pe09]} J. Peterson: {\it Examples of group actions which are virtually W*-superrigid}, 
 preprint arXiv:1002.1745.
\item {[PT07]} J. Peterson, A Thom: {\it Group cocycles and the ring of affiliated operators}, Invent. Math. {\bf 185} (2011), no. 3, 561--592.
\item {[PP86]} M. Pimsner, S. Popa: {\it Entropy and index for subfactors},  Ann. Sci. \'Ecole Norm. Sup. (4)  {\bf 19}  (1986),  no. 1, 57�-106.
\item {[Po01]} S. Popa: {\it On a class of type II$_1$ factors with Betti numbers invariants}, Ann. of Math. {\bf 163} (2006), 809--889. 
\item {[Po03]} S. Popa: {\it Strong Rigidity of II$_1$ Factors Arising from Malleable Actions of w-Rigid Groups. I.}, Invent. Math. {\bf 165} (2006), 369--408.
\item {[Po04]} S Popa: {\it Some computations of 1--cohomology groups and construction of non--orbit--equivalent
actions}, J. Inst. Math. Jussieu {\bf 5} (2) (2006), 309�-332.
\item {[Po05]} S. Popa: {\it Cocycle and orbit equivalence superrigidity for malleable actions of $w$-rigid groups},  Invent. Math.  {\bf 170}  (2007),  no. 2, 243--295.
\item {[Po06a]} S. Popa: {\it On the superrigidity of malleable actions with spectral gap}, 
J. Amer. Math. Soc. {\bf 21} (2008), 981--1000. 
\item {[Po06b]} S. Popa: {\it On Ozawa's Property for Free Group Factors}, Int. Math. Res. Notices (2007) Vol. 2007,
article ID rnm036.
\item {[Po07]} S. Popa: {\it Deformation and rigidity for group 	
actions and von Neumann algebras},  International Congress of Mathematicians. 
Vol. I,  445--477, 
Eur. Math. Soc., Z$\ddot{\text{u}}$rich, 2007.
\item {[Po09]} S. Popa: {Some results and problems in W$^*$--rigidity}, available at

 http://www.math.ucla.edu/popa/tamu0809rev.pdf.
\item {[PV09]} S. Popa, S. Vaes: {\it Group measure space decomposition of II$_1$ factors and} 

{\it W$^*$--superrigidity}, Invent. Math. {\bf 182} (2010), no. 2, 371--417.
\item {[PV11]} S. Popa, S. Vaes: {\it Unique Cartan decomposition for II$_1$ factors arising from arbitrary actions of free groups}, preprint arXiv:1111.6951.
\item {[PV12]} S. Popa, S. Vaes: {\it Unique Cartan decomposition for II$_1$ factors arising from arbitrary actions of hyperbolic groups},  preprint arXiv:1201.2824.
\item {[Si55]} I.M. Singer: {\it Automorphisms of finite factors}, Amer. J. Math. {\bf 77} (1955), 117--133.
 \item {[Si10]} T. Sinclair: {\it Strong solidity of group factors from lattices in SO(n,1) and SU(n,1)}, J. Funct. Anal., Volume {\bf 260} (2011), no.11, 3209--3221.

\item {[Ta03]} M Takesaki: {\it Theory of operator algebras, III,} Encyclopaedia of Mathematical Sciences, 127. Operator Algebras and Non-commutative Geometry, 8. Springer-Verlag, Berlin, 2003. xxii+548 pp.
\item {[Va07]} S. Vaes: {\it Explicit computations of all finite index bimodules for a family of II$_1$ factors}, Ann. Sci. \'Ec. Norm. Sup\'er. (4) {\bf 41} (2008), no. 5, 743--788.
\item {[Va10a]} S. Vaes: {\it Rigidity for von Neumann algebras and their invariants}
In Proceedings of the ICM (Hyderabad, India, 2010), Vol. III, Hindustan Book Agency (2010),  1624--1650.
\item {[Va10b]} S Vaes: {\it One--cohomology and the uniqueness of the group
measure space decomposition of a II$_1$ factor}, preprint arXiv:1012.5377, to appear in Math. Ann.
\enddocument